\documentclass[11pt,a4paper,notitlepage,oneside]{amsart}
\usepackage{etoolbox}
\usepackage{microtype} 
\usepackage{amsmath}
\usepackage{amsthm}
\usepackage{braket}
\usepackage{pdfpages}
\usepackage[psamsfonts]{amssymb}
\usepackage[T1]{fontenc}
\usepackage[english]{babel}
\usepackage[utf8]{inputenc}
\usepackage{enumitem}
\usepackage{mathtools}
\usepackage{caption,subfig} 
\usepackage{rotating} 
\usepackage{scalerel}
\usepackage{relsize} 
\usepackage{cancel} 

\usepackage[autostyle,italian=guillemets]{csquotes}
\usepackage[style=alphabetic,sorting=nyt,firstinits=true,maxnames=5,maxalphanames=5,backend=bibtex]{biblatex}
\bibliography{biblio}
\allowdisplaybreaks 

\usepackage{pgf,tikz,pgfplots}
\usepackage{mathrsfs}
\usetikzlibrary{arrows}

\usepackage{fancyhdr}
\usepackage{tikz-cd}
\usepackage{amscd}
\usepackage[a4paper, bottom=3cm, left=3cm,right=3cm]{geometry}
\usepackage{hyperref}
\hypersetup{hidelinks} 

\usepackage{subfiles}
\usepackage{my_shortcuts}

\newcommand{\quickref}{black} 
\newcommand{\tocheck}{black} 
\newcommand{\correction}{black}
\newcommand{\newcor}{black}

\author{Jan Hendrik Bruinier}
\address{
Fachbereich Mathematik, Technische Universität Darmstadt, Schlossgartenstraße 7, D–
64289 Darmstadt, Germany.
}
\email{bruinier@mathematik.tu-darmstadt.de}
\author{Riccardo Zuffetti}
\address{
Fachbereich Mathematik, Technische Universität Darmstadt, Schlossgartenstraße 7, D–
64289 Darmstadt, Germany.
}
\email{zuffetti@mathematik.tu-darmstadt.de}

\subjclass[2020]{11F27, 14J10, 11F46, 14C25.}

\title{Lefschetz decompositions of Kudla--Millson theta functions}

\begin{document}
	\maketitle
	\begin{abstract}
	In the \textcolor{\newcor}{1980's} Kudla and Millson introduced a theta function in two variables.
	It behaves as a Siegel modular form with respect to the first variable, and is a closed differential form on an orthogonal Shimura variety~$X$ with respect to the other variable.
	We prove that the Lefschetz decomposition of the cohomology class of that theta function is also its modular decomposition in Eisenstein, Klingen and cuspidal parts.
	We also show that the image of the Kudla--Millson lift is contained in the primitive cohomology of~$X$.
	As an application, we deduce dimension formulas for the cohomology of~$X$ in low degree.
	Our results cover the cases of moduli spaces of compact hyperkähler manifolds.
	\end{abstract}
	
	\tableofcontents
	\section{Introduction}\label{sec;outline}
	Orthogonal Shimura varieties provide a large family of varieties of geometric interest.
	Among them there are modular and Shimura curves, Hilbert modular surfaces, Siegel~$3$-folds, and moduli spaces of compact hyperkähler manifolds.
	
	Every orthogonal Shimura variety~$X$ comes with a large supply of algebraic cycles, given by embedded orthogonal Shimura subvarieties of smaller dimension.
	By celebrated work of Kudla--Millson, 
    the cohomology classes of such special cycles are related to Siegel modular forms via the geometric theta correspondence, see e.g.~\cite{kudlamillson-intersection} \textcolor{\newcor}{for an overview in greater generality. This correspondence} can be realized by an integral kernel, the so-called \emph{Kudla--Millson theta function}~$\KMtheta{\genus}$ of genus~$g$.
	It is a (non-holomorphic) Siegel modular form of genus~$g$ taking values in the space of closed differential~$(g,g)$-forms on~$X$.
	Its cohomology class~$[\KMtheta{\genus}]$ is a holomorphic Siegel modular form whose Fourier coefficients are the fundamental classes of codimension~$g$ special cycles on~$X$, generalizing the generating series of intersection numbers of Hirzebruch--Zagier divisors on Hilbert modular surfaces studied by Hirzebruch and Zagier in~\cite{hirzebruch-zagier}.

	Although~$X$ may be non-compact, the cohomology group~$H^{2\genus}(X)$ admits a \emph{Lefschetz decomposition} for small~$g$.
	Also the space of Siegel modular forms of large weight admits a natural decomposition, the \emph{modular decomposition}.
	It splits as a direct sum of subspaces given by cusp forms, Siegel Eisenstein series, and Klingen Eisenstein series.
	
	It is the goal of the present paper to provide the Lefschetz decomposition of the cohomology class~$[\KMtheta{\genus}]$, and to show that it is at the same time the modular decomposition of~$[\KMtheta{\genus}]$ as a Siegel modular form.
	More specifically, we prove that each term of the Lefschetz decomposition of~$[\KMtheta{\genus}]$ is given by the Klingen Eisenstein series arising from the primitive part of a Kudla--Millson theta function of intermediate genus.
	We also illustrate several geometric applications on orthogonal Shimura varieties.
	
%
%
	
	\subsection{Lefschetz and Modular decompositions}	
	To illustrate the main results of the present paper we need to introduce some notation.
	Let~$L$ be an even lattice of signature~$(n,2)$.
	To simplify the exposition, \textcolor{\newcor}{throughout} this introduction we assume that~$L$ is \emph{unimodular}.
	This condition will be dropped in the main body of the article.
	
	Let~$X$ be a (connected) orthogonal Shimura variety arising from~$L$, i.e.\ the quotient~${X=\Gamma\backslash D}$ for some finite index subgroup~$\Gamma\subseteq \OO^+(L)$, where~$D$ is the Hermitian symmetric domain of type IV attached to~$\OO(L)$ and~$\OO^+(L)$ is the group of isometries of the lattice~$L$ with positive spinor norm.
	Then~$X$ is a quasi-projective variety of dimension~$n$, which may have finite quotient singularities.
	If~$\Gamma$ is torsion free, then~$X$ is smooth.
	We denote by~$\omega$ the Kähler form on~$X$ arising as the first Chern form of the Hodge bundle.
	
	Let~$H^{2\genus}(X)$, resp.~$\sqH^{2\genus}(X)$, be the de Rham cohomology group, resp.~$L^2$-cohomology group, of~$X$ with complex coefficients.
	If~$\genus<(n-1)/2$, then the natural homomorphism
	\be\label{eq;introhomiso}
	\sqH^{2\genus}(X)\longrightarrow H^{2\genus}(X)
	\ee
	is an isomorphism preserving the Hodge structures.
	The \emph{Lefschetz decomposition} of~$\sqH^{\genus,\genus}(X)$ induces a Lefschetz decomposition on~$H^{\genus,\genus}(X)$ of the form
	\be
	H^{\genus,\genus}(X)\label{eq;lefdecintrodisp}
	=
	\CC[\omega]^\genus
	\oplus
	\lefschetz^{\genus-1} \Hprim^{1,1}(X)
	\oplus
	\cdots
	\oplus
	\lefschetz \Hprim^{\genus-1,\genus-1}(X)
	\oplus
	\Hprim^{\genus,\genus}(X).
	\ee
	Here we denoted by~$\Hprim^{*}(X)$ the primitive cohomology groups of~$X$ and by~$\lefschetz$ the Lefschetz operator given by wedging with~$\omega$; see Section~\ref{sec;digressdiffgeo} for further information.
	
	Let~$k=1+n/2$.
	From the well-known classification of unimodular lattices, it is easy to see that~$k$ is an \emph{even} integer.
	Let~$M^k_\genus$ be the space of weight~$k$ and genus~$\genus$ Siegel modular forms.
	The subspace of cusp forms of~$M^k_\genus$ is denoted by~$S^k_\genus$.
	
	From now on, we assume that~$k>2\genus$.
	The first classical example\textcolor{\newcor}{s} of element\textcolor{\newcor}{s} \textcolor{\newcor}{in}~$M^k_\genus$ \textcolor{\newcor}{are} the Siegel Eisenstein series~$E^k_\genus$. 
	These series have been generalized by Klingen~\cite{klingen-defKES} to the nowadays called Klingen Eisenstein series; see Section~\ref{sec;KMtheta} for details.
	They arise from genus~$r$ cusp forms with~$r\leq\genus$.
	We denote by~$E^k_{\genus,r}(f)$ the Klingen Eisenstein series of weight~$k$ and genus~$g$ arising from~$f\in S^k_r$.
	The subspace of such Klingen Eisenstein series is denoted by~$M^k_{\genus,r}$.
	We remark that in the degenerate case of~$g=r$, the series~$E^k_{\genus,\genus}(f)$ equals~$f$, and that if~$\genus=0$, then~$E^k_{\genus,0}(1)=E^k_\genus$.
	
	The space of modular forms~$M^k_\genus$ admits a decomposition
	\be\label{eq;moddecintro}
	M^k_\genus
	= \CC E^k_\genus
	\oplus
	M^k_{\genus,1}
	\oplus
	\dots
	\oplus
	M^k_{\genus,\genus-1}
	\oplus
	S^k_\genus,
	\ee
	which we call the \emph{modular decomposition} of~$M^k_\genus$.
	
	\subsection{The main result}
	The cohomology class~$[\KMtheta{\genus}]$ of the Kudla--Millson theta function may be regarded as an element of~$M^k_\genus\otimes H^{\genus,\genus}(X)$.
	It is then a natural problem to find its Lefschetz decomposition with respect to~\eqref{eq;lefdecintrodisp}, as well as its modular decomposition with respect to~\eqref{eq;moddecintro}.
	
	The main result of the present paper, stated below, answers these questions.
	For~$0\leq r \leq \genus$, we denote by~$\klin_{\genus,r}^k\colon M^k_r\to M^k_\genus$ the linear map defined with respect to the modular decomposition of~$M^k_r$ as
	\[
	\klin_{\genus,r}^k\colon
	E^k_{r,i}(f_i)
	\longmapsto
	E^k_{\genus,i}(f_i),
	\qquad
	\text{\textcolor{\newcor}{for every~$i\leq r$ and} every~$f_i\in S^k_i$.}
	\]
	\begin{thm}\label{thm;mainresintro}
	Let~$g<(n+2)/4$ be a positive integer.
	The \emph{primitive part} of~$[\KMtheta{\genus}]$ is
	\be\label{eq;primpartintro}
	[\KMtheta{\genus}]_0
	=
	[\KMtheta{\genus}]
	+
	\lefschetz \circ \klin^k_{\genus,\genus-1}\big([\KMtheta{\genus-1}]\big)
	\ee
	and lies in~$S^k_\genus\otimes\Hprim^{\genus,\genus}(X)$.
	In particular, it is a Siegel cusp form with respect to the symplectic variable.
	The \emph{Lefschetz decomposition} of~$[\KMtheta{\genus}]$ equals the \emph{modular decomposition} of~$[\KMtheta{\genus}]$, and is given by
	\be\label{eq;lefdecintro}
	[\KMtheta{\genus}]
	=
	\sum_{r=0}^\genus
	(-\lefschetz)^{\genus-r}
	\circ
	E^k_{\genus,r}\big([\KMtheta{r}]_0\big).
	\ee
	Here the summand of index~$r$ lies in~$M^k_{\genus,r}\otimes\lefschetz^{\genus-r}\Hprim^{r,r}(X)$.
	\end{thm}
    \textcolor{\newcor}{For an explanation of the obstacles to weakening the hypothesis~$g<(n+2)/4$ of Theorem~\ref{thm;mainresintro}, see Section~\ref{sec:intro:constraints}.}
    
	
	In this paper we prove a more general version of Theorem~\ref{thm;mainresintro}, in which~$L$ is not assumed to be unimodular; see Theorem~\ref{thm;lefschetzdecgengenus}.
	To do so, it is necessary to deal with vector-valued Siegel modular forms associated to the Weil representation, which are introduced in Section~\ref{sec;KMtheta}.
	
	We now provide an idea of the proof of Theorem~\ref{thm;mainresintro}.
	The first step is to compute the modular decomposition of~$[\KMtheta{\genus}]$ by induction on~$\genus$.
	We show that the image of~$[\KMtheta{\genus}]$ under the Siegel~$\Phi$ operator equals (up to a sign) the image of~$[\KMtheta{\genus-1}]$ under the Lefschetz operator.
	This is used to prove that the right-hand side of~\eqref{eq;primpartintro} is the cuspidal part of~$[\KMtheta{\genus}]$.
	
	The technical heart of the proof is to show that each term of the modular decomposition of~$[\KMtheta{\genus}]$ lies in a different term of the Lefschetz decomposition of~$H^{\genus,\genus}(X)$.
	To achieve this goal, it is sufficient to prove that the image of the \emph{Kudla--Millson lift}
	\[
	\KMlift{\genus}\colon S^k_\genus\longrightarrow H^{\genus,\genus}(X),
	\qquad
	f\longmapsto\pet{f}{[\KMtheta{\genus}]},
	\]
	is contained in the primitive cohomology~$\Hprim^{\genus,\genus}(X)$.
	Here~$\petempty$ denotes the Petersson inner product on~$S^k_\genus$.
	Theorem~\ref{thm;mainresintro} may then be deduced from the following result of independent interest.
	\begin{thm}\label{thm;imageKMliftintro}
	Let~$\genus<(n+2)/4$.
	The image of the Kudla--Millson lift~$\KMlift{\genus}$ is a subspace of the primitive cohomology group~$\Hprim^{\genus,\genus}(X)$.
	If the codimension~$\genus$ special cycles on~$X$ arising from~$L$ span~$H^{\genus,\genus}(X)$, then the image of~$\KMlift{\genus}$ is the whole~$\Hprim^{\genus,\genus}(X)$.
	\end{thm}
	We remark that Theorem~\ref{thm;imageKMliftintro} corrects \textcolor{\quickref}{the} statement on the image of the Kudla--Millson lift made in~\cite[Theorem, p.4]{kudlamillson-tubes}; see Section~\ref{sec;imageKMlift} for details.
    
    \textcolor{\newcor}{Using the results of Kudla--Millson \cite{kudlamillson-intersection} and the Noether--Lefschetz Conjecture, proved  in~\cite{noetherlefconj},  Theorem~\ref{thm;imageKMliftintro} can be reduced to showing} that
	\be\label{eq;vanintegintro}
	\int_X \KMlift{\genus}(f)\wedge\KMtheta{\genus-1}(\tau)\wedge \omega^{n-(2\genus-1)}=0
	\ee
	for all~$f\in S^k_\genus$ and all~$\tau\in \HH_{\genus-1}$.
	By means of the doubling method, we prove  that~\eqref{eq;vanintegintro} follows from the vanishing of Garrett's map
	\[
	S^k_\genus\longrightarrow M^k_{\genus-1},\qquad
	f\longmapsto\bigpet{f}{E^k_{2\genus-1}\big(\begin{smallmatrix}
	\tau & 0 \\
	0 & \ast
	\end{smallmatrix}\big)},
	\]
	where we denote by~$\ast$ the variable of integration in~$\HH_\genus$.
	
	In Section~\ref{sec;adelicspcy} of this paper we also provide criteria on~$L$ so that the codimension~$g$ special cycles arising from~$L$ span~$H^{\genus,\genus}(X)$.
	Note that this does not immediately follow from~\cite{noetherlefconj}, due to a small inaccuracy of~\cite[Remark~$2.4$]{noetherlefconj}; see Remark~\ref{rem;BLMMhasimprecision} for details.
	We correct it with a careful study of the irreducible components of adelic special cycles.
	
	\subsection{Applications of geometric interest}
	Immediate consequences of our main results are the following bounds for the dimension of cohomology groups of orthogonal Shimura varieties.
	\begin{cor}\label{cor;ifinjthenformcohointro}
	Let~$\genus<(n+2)/4$.
	Then
	\[
	\dim\lefschetz^{\genus-r} H_0^{r,r}(X)\geq\dim \KMlift{r}\big(S^k_{r}\big)\qquad\text{for every~$0\leq r\leq\genus$}.
	\]
	If the Kudla--Millson lift~$\KMlift{r}$ is injective for every~$r\leq \genus$, then
	\[
	\dim\lefschetz^{\genus-r} H_0^{r,r}(X)\geq\dim S^k_{r}
	\qquad\text{and}\qquad
	\dim H^{\genus,\genus}(X)\geq\dim M^k_\genus.
	\]
	\end{cor}
	
	The injectivity of~$\KMlift{\genus}$ in the case of genus~$\genus=1$ has been proved in~\cite{bruinierbook}, \cite{bruinierfunke}, \cite{zuffetti-unfolding}.
%
	Recently, Kiefer and the second author proved in~\cite{kieferzuffetti} the injectivity of~$\KMlift{\genus}$ in the case of genus~$\genus=2$ under the assumption that~$n>9$.
	The papers cited above provide criteria on the injectivity of the lift that work even without assuming~$L$ to be unimodular.
	
	As an interesting example, we consider the moduli spaces~$X_d$ of degree~$2d$ quasi-polarized~K3 surfaces.
	They can be viewed as orthogonal Shimura varieties arising from even lattices of the form
	\[
	L_d=\langle 2d\rangle \oplus U\oplus U\oplus E_8\oplus E_8,
	\]
	where~$U$ is the hyperbolic lattice of signature~$(1,1)$ and~\textcolor{\newcor}{$E_8$ is the root lattice of the~$E_8$ root system}.
	The lattices~$L_d$ are not unimodular, but satisfy the known injectivity criteria for~$\KMlift{1}$ and~$\KMlift{2}$.
	Since the spaces~$H^{1,1}(X_d)$ and~$H^{2,2}(X_d)$ are spanned by classes of special cycles arising from~$L_d$, we may deduce the following result from Theorem~\ref{thm;imageKMliftintro} and Corollary~\ref{cor;ifinjthenformcohointro}.
	Let~$S^k_{\genus,\lattice}$ be the space of weight~$k$ and genus~$\genus$ vector-valued cusp forms with respect to the Weil representation arising from~$\lattice$.
	\begin{cor}\label{cor;cohomodspintro}
	Let~$X_d$ be the moduli space of degree~$2d$ quasi-polarized K3 surfaces.
	Then
	\[
	\dim H^2(X_d)= 1 + \dim S^k_{1,L_d}
	\qquad\text{and}\qquad
	\dim H^4(X_d)=1+\dim S^k_{1,L_d} + \dim S^k_{2,L_d}.
	\]
	\end{cor}
	The dimension formula for~$H^2(X_d)$ above was proved in~\cite{bruinier-converse} and~\cite{noetherlefconj}.
	Our article provides a different proof, which has the advantage of working also for cohomology groups of higher degree.

\subsection{Constraints on the genus}\label{sec:intro:constraints}
\textcolor{\newcor}{The main results illustrated above hold if the genus~$g$ of the theta series satisfies~$g<(n+2)/2$.
There are three main obstacles to weakening this assumption.}

\textcolor{\newcor}{The first is the required  convergence of the Klingen Eisenstein series appearing in~\eqref{eq;lefdecintro}, which holds only if~$k>2\genus$.
	Since~$k=1+n/2$, the latter bound is equivalent to~$g<(n+2)/4$.
	 To bypass this problem, one could study the analytic continuation of Klingen Eisenstein series with respect to spectral parameters, as in \cite{weissauer-88}. This leads to the subtle issue of proving holomorphicity in the symplectic variable of their continuations.}

\textcolor{\newcor}{The second obstacle is the Noether--Lefschetz Conjecture.
This is proved in~\cite{noetherlefconj} under the assumption that~$g<(n+1)/3$.
The fact that the codimension~$g$ special cycles span~$H^{g,g}(X)$ is used to reduce the proof of Theorem~\ref{thm;imageKMliftintro} to~\eqref{eq;vanintegintro}. As pointed out in Remark~2.8 of~\cite{noetherlefconj} it is expected that this bound is sharp. For instance, already in the case of Hilbert modular surfaces (where $n=2$) the middle cohomology is in general not only generated by Hirzebruch--Zagier divisors. It was shown by Klingenberg and Murty--Ramakrishnan that there are additional Hodge classes corresponding to certain Hilbert cusp forms of weight $2$ with complex multiplication. 
 }

\textcolor{\newcor}{The third is a comparison of the de Rham and~$L^2$-cohomology groups of~$X$.
Under the slightly weaker bound~$g<(n-1)/2$, the homomorphism~\eqref{eq;introhomiso} is an isomorphism and induces a Lefschetz decomposition of~$H^{2g}(X)$ and hence of~$[\KMtheta{g}]$.
It is unclear whether, for higher~$g$, some decomposition of~$H^{2g}(X)$ induces a decomposition of~$[\KMtheta{g}]$ with automorphic meaning.
To avoid the comparison between $L^2$ and de Rham cohomology, one could replace~$X$ with a compact orthogonal Shimura variety defined over a totally real number field of degree greater than~$1$, so that the Lefschetz decomposition would be available on the full cohomology.
However, the two other constrains mentioned above would still remain.
}

	\subsection*{Acknowledgments}
    \textcolor{\newcor}{We thank the referees for their careful reading of our paper and for helpful comments.}
	We would like to thank Siegfried Böcherer, Soumya Das, Jens Funke, \textcolor{\newcor}{Luis García}, Paul Kiefer, Zhiyuan Li, \textcolor{\newcor}{Niklas Ludwig,} Martin Möller, Riccardo Salvati-Manni, Oliver Stein and Brandon Williams for useful conversations on the topics of this paper.
	The authors acknowledge support by Deutsche Forschungsgemeinschaft (DFG, German Research Foundation) through the Collaborative Research Centre TRR 326 \textit{Geometry and Arithmetic of Uniformized Structures}, project number 444845124.
	\section{Vector-valued Siegel modular forms}\label{sec;KMtheta}
	In this section we gather all the properties of vector-valued Siegel modular forms with respect to the Weil representation relevant for the purposes of the paper.
	Although the following theory of vector-valued modular forms can be extended to lattices of general signature, we prefer to work in signature~$(n,2)$, since this is the only case needed in the rest of the present article.
	For the same reason, we only consider weights~$k\in\frac{1}{2}\ZZ$ such that~$2k\equiv n-2$ mod~$4$.
	
	\subsection{The Weil representation}
	\textcolor{\newcor}{We denote by~$M^t$ the transpose of any matrix~$M$.
    If~$M$ is invertible with entries in a commutative ring, then we write~$M^{-t}\coloneqq(M^t)^{-1}=(M^{-1})^t$.}
    
	The Siegel upper half-space of genus~$\gengenus$ is defined as
	\[
	\HH_{\gengenus}=\{\tau=x+iy\in\CC^{\gengenus\times \gengenus} : \text{$\tau=\tau^t$ and $y>0$}\}.
	\]
	Here we denote by~$x$ and~$y$ respectively the real and imaginary part of~$\tau$.
	
	If~$z\in\CC$, we write~$e(z)=\exp(2\pi i z)$ and denote by~$z^{1/2}=\sqrt{z}$ the principal branch of the square root, so that~$\arg(\sqrt{z})\in(-\pi/2,\pi/2]$.
	For every integer~$m$ we define~$z^{m/2}=(z^{1/2})^m$.
	
	We write~$\Mp_{2\gengenus}(\RR)$ for the metaplectic cover of the symplectic group~$\Sp_{2\gengenus}(\RR)$.
	This is realized by choosing a square root of the map~$\tau\mapsto \det(C\tau+D)$ for every~$\big(\begin{smallmatrix}
	A & B\\
	C & D
	\end{smallmatrix}\big)\in\Sp_{2\gengenus}(\RR)$.
	The group law of~$\Mp_{2\gengenus}(\RR)$ is
	\be\label{eq;grouplawmetgr}
	\big(M_1,\phi_1(\tau)\big)\cdot \big(M_2,\phi_2(\tau)\big)
	=
	\big(M_1M_2 , \phi_1\big(M_2\cdot\tau\big)\phi_2(\tau)\big),
	\ee
	where the action of~$\Mp_{2\gengenus}(\RR)$ on~$\HH_\gengenus$ is the one induced by the standard action of~$\Sp_{2\gengenus}(\RR)$.
	
	The standard generators of the group~$\Mp_{2\gengenus}(\ZZ)$ are the elements of the form
	\be\label{eq;stdgenofmetcov}
	S=\left(\left(\begin{smallmatrix}
	0 & -I_\gengenus\\
	I_\gengenus & 0
	\end{smallmatrix}\right),\sqrt{\det\tau}\right)
	\qquad\text{and}\qquad
	T_B=\left(\left(\begin{smallmatrix}
	I_\gengenus & B \\
	0 & I_\gengenus
	\end{smallmatrix}\right),1\right)\quad\text{where~$B\in\Sym_\genus(\ZZ)$.}
	\ee
	For future purposes we introduce here also the elements of~$\Mp_{2\gengenus}(\ZZ)$ given by
	\be\label{eq;defRgenmetapl}
	R_A=\left(\left(\begin{smallmatrix}
	A & 0\\
	0 & A^{-t}
	\end{smallmatrix}\right), \det A^{1/2}\right)\qquad\text{for every~$A\in\GL_\gengenus(\ZZ)$.}
	\ee
	
	Let~$(\lattice,q)$ be an even indefinite lattice of signature~$(n,2)$, and let~$L'$ be its dual lattice.
	We denote by~$(\cdot{,}\cdot)$ the bilinear form associated to~$q$ such that~$q(\lambda)=(\lambda,\lambda)/2$ for every~$\lambda\in\lattice$.
	If~$\lambda\in \lattice^\genus$, we denote by~$q(\lambda)$ the~$\genus\times\genus$ moment matrix of the tuple~$\lambda$.
	The discriminant group of~$\lattice$ is the quotient~$D_\lattice=\lattice'/\lattice$.
	The form~$q$ induces a~$\QQ/\ZZ$-valued moment form on~$D_\lattice^\genus$ which we still denote by~$q$.
	
	We endow the group algebra~$\CC[D_\lattice]$ with the Hermitian scalar product~$\langle\cdot{,}\cdot\rangle$ defined as
	\be\label{eq;stinprodgroupalg}
	\Big\langle
	\sum_{\discel\in D_\lattice}\lambda_\discel \mathfrak{e}_\discel
	,
	\sum_{\discel\in D_\lattice}\mu_\discel \mathfrak{e}_\discel
	\Big\rangle
	=
	\sum_{\discel\in D_\lattice}\lambda_\discel\overline{\mu_\discel},
	\ee
	where~$(\mathfrak{e}_\discel)_\discel$ is the standard basis of~$\CC[D_\lattice]$.
	We denote with the same symbol also the extension of that scalar product to~$\CC[D_\lattice^\gengenus]$ which makes different copies of~$D_\lattice$ orthogonal to each others. 
	
	\begin{defi}
	The \emph{Weil representation} of genus~$\gengenus$ associated to the lattice~$\lattice$ is the unitary representation
	\[
	\weil{\lattice}{\gengenus}\colon \Mp_{2\gengenus}(\ZZ) \longrightarrow \Aut(\CC[D_\lattice^\gengenus])
	\]
	defined by
	\begin{align*}
	\weil{\lattice}{\gengenus}(T_B)\mathfrak{e}_\discel &=
	e\big(
	\trace (q(\discel)B)
	\big)\mathfrak{e}_\discel\qquad\text{for every~$B\in\Sym_\gengenus(\ZZ)$,}
	\\
	\weil{\lattice}{\gengenus}(S)\mathfrak{e}_\discel &=
	\frac{\sqrt{i}^{\gengenus(2-n)}}{|D_\lattice|^{\gengenus/2}}
	\sum_{\disceltwo\in D_\lattice^\gengenus} e\big(-\tr ( \disceltwo,\discel)\big)\mathfrak{e}_\disceltwo,
	\\
	\weil{\lattice}{\gengenus}(R_A)\mathfrak{e}_\discel
	&=
	(\det A)^{(2-n)/2}\mathfrak{e}_{\discel A^{-1}}\qquad
	\text{for every~$A\in\GL_\gengenus(\ZZ)$},
	\end{align*}
	see e.g.~\cite{borcherds-grassmannians} and~\cite{zh;phd}.
	\end{defi}
	
	We remark that
	\be\label{eq;howto-phi}
	\weil{\lattice}{\gengenus}\big((M,-\phi)\big)\mathfrak{e}_\discel
	=
	(-1)^{2-n}\weil{\lattice}{\gengenus}\big( (M,\phi)\big)\mathfrak{e}_\discel
	\qquad\text{for every~$(M,\phi)\in\Mp_{2\gengenus}(\ZZ)$.}
	\ee
	This follows from the fact that~$(M,-\phi)=(M,\phi)\cdot R_A^2$, where~$R_A$ is chosen as in~\eqref{eq;defRgenmetapl} with~$A=\big(
	\begin{smallmatrix}
	-1 & \\
	 & I_{\gengenus-1}
	\end{smallmatrix}
	\big)$.
	
	Weil representations of different genera associated to the same lattice~$L$ are \textcolor{\newcor}{closely connected}.
	To clarify this, we need to introduce another piece of notation.
	We identify~$\CC[D_\lattice^r]\otimes\CC[D_\lattice^{\gengenus-r}]$ with~$\CC[D_\lattice^\gengenus]$ under the isomorphism
	\be\label{eq;identgroupalg}
	\CC[D_\lattice^r]\otimes\CC[D_\lattice^{\gengenus-r}]\longrightarrow\CC[D_\lattice^\gengenus],\qquad \mathfrak{e}_\discel\otimes\mathfrak{e}_\disceltwo\longmapsto\mathfrak{e}_{(\discel,\disceltwo)}.
	\ee
	Let~$\iota$ be the inclusion of metaplectic groups given as
	\ba\label{eq;inclprodsimpl}
	\iota\colon\Mp_{2r}(\ZZ)\times\Mp_{2(\gengenus - r)}(\ZZ)&\longrightarrow\Mp_{2\gengenus}(\ZZ),
	\\
	\Big(\big(\big(\begin{smallmatrix}
	A & B\\
	C & D
	\end{smallmatrix}
	\big),\phi\big),\big(\big(
	\begin{smallmatrix}
	A' & B'\\
	C' & D'
	\end{smallmatrix}\big),\phi'\big)\Big)
	&\longmapsto
	\left(\left(\begin{smallmatrix}
	A & & B & \\
	 & A' & & B' \\
	C & & D & \\
	 & C' & & D'
	\end{smallmatrix}\right),\tilde{\phi}\right),
	\ea
	where the function~$\tilde{\phi}(\tilde{\tau})$, with~$\tilde{\tau}\in\HH_\gengenus$, is~$\pm\det\big(\big(\begin{smallmatrix}
	C & \\
	 & C'
	\end{smallmatrix}\big)\tilde{\tau} + \big(\begin{smallmatrix}
	D & \\
	 & D'
	\end{smallmatrix}\big)\big)^{1/2}$ with sign chosen such that~$\tilde{\phi}\big(\begin{smallmatrix}
	\tau & 0 \\
	0 & \tau'
	\end{smallmatrix}\big)= \phi(\tau)\cdot\phi'(\tau')$.
	
	Weil representations of different genera attached to the same lattice behave well with respect to the inclusion~$\iota$.
	This is illustrated in the following result, whose proof is left to the reader and may be deduced with the same idea of~\cite[Lemma~$3.2$]{kieferzuffetti}.
	\begin{lemma}\label{lemma;proofweilreponsubmeta}
	Let~$0\le r\le g$, and let~$\iota$ be the inclusion as in~\eqref{eq;inclprodsimpl}.
	If~$\gamma\in\Mp_{2r}(\ZZ)$ and~$\gamma'\in\Mp_{2(\gengenus-r)}(\ZZ)$, then
	\be\label{eq;inlcfmet2}
	\weil{\lattice}{\gengenus}\big(\iota(\gamma,\gamma')\big)(\mathfrak{e}_\discel\otimes\mathfrak{e}_\disceltwo)
	=
	\big(\weil{\lattice}{r}(\gamma)\mathfrak{e}_\discel\big)\otimes
	\big(\weil{\lattice}{\gengenus-r}(\gamma')\mathfrak{e}_\disceltwo\big)
	\ee
	for every~$\discel\in D_\lattice^r$ and~$\disceltwo\in D_\lattice^{\gengenus-r}$.
	\end{lemma}

	\subsection{Vector-valued modular forms}\label{sec;backgrvvmodforms}
	Let~$k\in\frac{1}{2}\ZZ$ be such that~$2k\equiv n-2$ mod~$4$.
	A Siegel modular form of weight~$k$ and genus~$\gengenus$ with respect to the Weil representation~$\weil{\lattice}{\gengenus}$ is a holomorphic function~$f\colon\HH_\gengenus\to\CC[D_\lattice^\gengenus]$ such that
	\[
	f(\gamma\cdot\tau)=\phi(\tau)^{2k} \weil{\lattice}{\gengenus}(\gamma) f(\tau)
	\qquad
	\text{for every~$\gamma=\big(M , \phi(\tau)\big)\in\Mp_{2\gengenus}(\ZZ)$},
	\]
	with the usual requirement that~$f$ is holomorphic at the cusp~$\infty$ for the case of genus~$1$.
	We write~$f=\sum_{\discel\in D_\lattice^\gengenus}f_\discel\mathfrak{e}_\discel$, denoting by~$f_\discel$ the~$\discel$-component of~$f$.
	Such components are scalar-valued Siegel modular forms of level~$N$, where~$N$ is the level of~$\lattice$, i.e.\ the smallest positive integer such that~$N q(\lambda)\in\ZZ$ for every~$\lambda\in \lattice'$.
	
	Let~$\halfint{\gengenus}$
	be the set of symmetric half-integral~$\gengenus\times\gengenus$-matrices.
	The invariance of~$f$ with respect to the translations on~$\HH_\gengenus$ induced by~$\Mp_{2\gengenus}(\ZZ)$ implies that~$f$ admits a Fourier expansion of the form
	\be\label{eq;Fexpfirstappvvc}
	f(\tau)=\sum_{\discel\in D_\lattice^\gengenus}
	\sum_{\substack{T\in q(\discel) + \halfint{\gengenus} \\ T\ge 0}}
	c_T(f_\discel) q^T \mathfrak{e}_\discel,
	\ee
	where~$q^T=e \big(\tr (T\tau) \big)$.
	If~$\genus=1$, then the fact that the series~\eqref{eq;Fexpfirstappvvc} is indexed over non-negative rational numbers is implied by \textcolor{\newcor}{requiring} that elliptic modular forms are holomorphic at~$\infty$.
	For larger genus, one can show that the Fourier expansion is indexed only over positive semidefinite matrices by the Koecher principle for Siegel modular forms with level.
	
	We denote by
	\[
	c_T(f)=\sum_{\discel\in D_\lattice^\gengenus} c_T(f_\discel) \mathfrak{e}_\discel
	\]
	the vector in~$\CC[D_\lattice^\gengenus]$ gathering all Fourier coefficients of index~$T$ of~$f$.
	The Fourier coefficients of Siegel modular forms with respect to~$\weil{\lattice}{\gengenus}$ satisfy symmetric relations with respect to the action~$T\mapsto A^t T A$ induced by matrices~$A\in\GL_\gengenus(\ZZ)$.
	More precisely, we have
	\be\label{eq;symmofcoef}
	c_T(f) = \phi^{2k} \weil{\lattice}{\gengenus}\big(\big(\begin{smallmatrix}
	A & \\
	& A^{-t}
	\end{smallmatrix}\big),\phi\big) c_{A^t T A}(f)\qquad \text{for all~$A\in\GL_{\gengenus}(\ZZ),$}
	\ee
	where~$\phi$ is any square-root of~$\det A$.
	It is easy to see that~\eqref{eq;symmofcoef} implies that
	\be\label{eq;symmofcoef2}
	c_T(f_\discel)=c_{A^t T A}(f_{\discel A})
	\qquad
	\text{for all~$\discel\in D_\lattice^\genus$}.
	\ee
	
	If all Fourier coefficients~$c_T(f)$ of~$f$ with respect to singular matrices~$T$ vanish, then~$f$ is said to be a \emph{Siegel cusp form}.	
	We denote by~$M^k_{\gengenus,\lattice}$ and~$S^k_{\gengenus,\lattice}$ respectively the space of Siegel modular forms and cusp forms of weight~$k$ and genus~$\gengenus$ with respect to~$\weil{\lattice}{\gengenus}$, with the convention that~$S^k_{0,\lattice}=M^
	k_{0,\lattice}=\CC$ as in~\cite{klingen}.	

	\subsection{Relations with sublattices}\label{sec;relwithsublat}
	Let~$M\subset \lattice$ be a sublattice of finite index.
	There are two natural maps that relate the spaces~$M^k_{\genus,\lattice}$ and~$M^k_{\genus,M}$, as illustrated in the following result.
	This may be considered as a generalization to higher genus of~\cite[Section~$4$]{scheit} and~\cite[Lemma~$3.1$]{bruinier-yang09}.
	The natural surjective map~$\lattice'/M\to\lattice'/\lattice$, $\mu\mapsto\overline{\mu}$, induces a map~$(\lattice'/M)^\genus\to(\lattice'/\lattice)^\genus$ for every~$\genus$.
	\begin{lemma}\label{lemma;genofBYandS}
	There are two natural maps
	\[
	\res_{\genus,\lattice/M}\colon M^k_{\genus,\lattice}\longrightarrow M^k_{\genus,M},
	\qquad
	f\longmapsto f_M
	\]
	and
	\[
	\tr_{\genus,\lattice/M}\colon M^k_{\genus,M}\longrightarrow M^k_{\genus,\lattice},
	\qquad
	h\longmapsto h^\lattice
	\]
	such that~$\langle f, h^\lattice\rangle=\langle f_M,h\rangle$ for every~$f\in M^k_{\genus,\lattice}$ and~$h\in M^k_{\genus,M}$.
	They are given as
	\[
	\big(f_M\big)_\mu
	=
	\begin{cases}
	f_{\overline{\mu}}, & \text{if~$\mu\in (\lattice'/M)^\genus$,}
	\\
	0, & \text{otherwise,}
	\end{cases}
	\qquad
	\text{for every~$f\in M^k_{\genus,\lattice}$ and~$\mu\in D_M^\genus$}
	\]
	and
	\[
	\big(h^\lattice\big)_{\overline{\mu}}
	=
	\sum_{\discel\in (\lattice/M)^\genus}
	h_{\discel + \mu}
	\qquad
	\text{for every~$h\in M^k_{\genus,M}$,}
	\]
	where~$\mu$ is any fixed preimage of~$\overline{\mu}$ in~$(\lattice'/M)^\genus$.
	\end{lemma}
	\begin{proof}
	It is analogous to~\cite[Theorem~$4.1$]{scheit} and~\cite[Lemma~$3.1$]{bruinier-yang09}.
	\end{proof}

    \section{Klingen Eisenstein series and the doubling method}
    \textcolor{\quickref}{In this section we generalize the theory of Klingen Eisenstein series and Poincaré series to the case of vector-valued modular forms with respect to the Weil representation.
    Eventually, we illustrate the doubling method introduced by Garrett~\cite{garrett} and Böcherer~\cite{doubling} in this more general setting.
    This is a key ingredient to prove the modular decomposition of the Kudla--Millson theta function.}

	\subsection{Klingen Eisenstein series}\label{sec;KEseries}
	Here we introduce vector-valued Klingen Eisenstein series with respect to the Weil representation.
	In the case of unimodular lattices, these series are the classical scalar-valued Klingen Eisenstein series of~\cite[Chapter~II, Section~$5$]{klingen}.
	
	Let~$r\le\gengenus$ be a non-negative integer.
	We denote by~$\klipar{\gengenus}{r}$ and~$\klimet{\gengenus}{r}$ respectively the Klingen parabolic subgroup
	\[
	\klipar{\gengenus}{r}\coloneqq\left\{
	\gamma=\left(\begin{smallmatrix}
	A & B\\
	C & D
	\end{smallmatrix}\right)\in\Sp_{2\gengenus}(\ZZ) :
	A=\left(\begin{smallmatrix}
	\ast & 0\\
	\ast & \ast
	\end{smallmatrix}\right),
	C=\left(\begin{smallmatrix}
	\ast & 0\\
	0 & 0
	\end{smallmatrix}\right),
	D=\left(\begin{smallmatrix}
	\ast & \ast\\
	0 & \ast
	\end{smallmatrix}\right)
	\right\}
	\]
	and its preimage in~$\Mp_{2\gengenus}(\ZZ)$, where the top-left corners of the block decompositions of~$A$, $C$ and~$D$ above are~$r\times r$ matrices.
	
    
	We consider the projections
	\ba
	\ast&\colon\HH_\gengenus \longrightarrow\HH_r,
	\quad
	&&
	\tau=\big(\begin{smallmatrix}
	\tau^*  & \ast\\
	\ast & \ast
	\end{smallmatrix}\big)\longmapsto \tau^*,
	\\
	\ast&\colon \klipar{\gengenus}{r}\longrightarrow \Sp_{2r}(\ZZ),\quad
	&& M=\big(\begin{smallmatrix}
	A & B \\
	C & D
	\end{smallmatrix}\big)
	\longmapsto
	M^*=\big(\begin{smallmatrix}
	A_1 & B_1\\
	C_1 & D_1
	\end{smallmatrix}\big),
	\\
	\ast & \colon\klimet{\gengenus}{r}\longrightarrow \Mp_{2r}(\ZZ),\quad
	&&
	\big(M,\phi(\tau)\big)\longmapsto \big(M^*,\phi^*(\tau^*)\big),
	\ea
	where~$\phi^*(\tau^*)$ is the square root function chosen such that~$\phi(\tau)
	={\det D_4} ^{1/2} \phi^*(\tau^*)$.
	Here $A_1,\dots,D_1$ are the upper-left~$r\times r$ blocks of~$A,\dots,D$, and~$D_4$ is the lower-right~$(\gengenus-r)\times(\gengenus-r)$ block of~$D$.
	We remark that if~$M\in\klipar{\gengenus}{r}$, then
	\bas
	\phi(\tau)^2
	&=
	\det\left(
	\left(\begin{smallmatrix}
	C_1 & 0 \\
	0 & 0
	\end{smallmatrix}\right)\left(\begin{smallmatrix}
	\tau^* & \ast \\
	\ast & \ast
	\end{smallmatrix}\right) + \big(\begin{smallmatrix}
	D_1 & D_2 \\
	0 & D_4
	\end{smallmatrix}\big)
	\right)
	= \det D_4 \cdot \det(C_1\tau^* + D_1)
	= \det D_4 \cdot \phi^*(\tau^*)^2
	,
	\eas
	hence~$\phi^*(\tau^*)$ is defined to be equal to~$\phi(\tau)$ up to a square root of~$\det D_4$, which we choose to be the principal one.
	Note also that~$(\gamma\cdot \tau)^* = \gamma^*\cdot\tau^*$ for every~$\gamma\in \klipar{\gengenus}{r}$ and~$\tau\in\HH_\gengenus$.
	If it is unclear on which space~$\mathfrak{e}_0$ lies, we write~$\zeroe{r}$ to remark that it lies in~$\CC[D_\lattice^r]$.
	\begin{defi}\label{def;KlingEisseries}
	Let~$\gengenus,r\in\ZZ_{\ge0}$ be such that~$r\le\gengenus$.
	The \emph{Klingen Eisenstein series} of genus~$\gengenus$ arising from the cusp form~$f\in S^k_{r,\lattice}$ is defined as
	\[
	\klinser{k}{\gengenus}{r}{\lattice}(\tau,f)
	=
	\sum_{\gamma\in\klimet{\gengenus}{r}\backslash \Mp_{2\gengenus}(\ZZ)}
	\phi(\tau)^{-2k}
	\weil{\lattice}{\gengenus}(\gamma)^{-1}
	\Big(f\big(\textcolor{\newcor}{(\gamma\cdot\tau)^*}\big)
	\otimes \zeroe{\gengenus-r}\Big).
	\]
	The \emph{Siegel Eisenstein series} of genus~$g$ and weight~$k$ with respect to~$\weil{\lattice}{\genus}$
	is defined as~$\eiszero{k}{\gengenus}{\lattice}(\tau)\coloneqq\klinser{k}{\gengenus}{0}{\lattice}(\tau,1)$.
	\end{defi}
	It is easy to see by using \textcolor{\quickref}{e.g.~\cite[Proposition~$4.10$~(i)]{namikawa}} that the series~$\klinser{k}{\gengenus}{r}{\lattice}(f)$ is well-defined.
	A standard argument shows that~$\eiszero{k}{\gengenus}{\lattice}$ is normalized such that its Fourier coefficients of index~$T=0$ are
	\[
	c_0\big((\eiszero{k}{\gengenus}{\lattice})_\discel\big)
	=\begin{cases}
	1, & \text{if $\discel=0\in D_\lattice^\gengenus$,}\\
	0, & \text{otherwise}.
	\end{cases}
	\]
	
	The space of Siegel Eisenstein series with respect to~$\weil{\lattice}{\gengenus}$ can be of dimension greater than~$1$; see e.g.~\cite[Section~$1.2.3$]{bruinierbook} for the case of genus~$1$.
	For the purposes of the present article, we do not need to introduce other Eisenstein series than~$\eiszero{k}{\gengenus}{\lattice}$, i.e.\ the one arising from the zero of~$D_\lattice^\gengenus$.

	Scalar-valued Eisenstein series have excellent convergence properties.
	These are satisfied also by the vector-valued Klingen Eisenstein series defined above, as illustrated in the following generalization of~\cite[p.~$67$, Theorem~$1$]{klingen}, \textcolor{\quickref}{whose proof is standard and left to the reader}.
	For~$c>0$, we define the \emph{vertical strip of height~$c$} as
	\[
	V_\gengenus(c)
	\coloneqq
	\{
	\tau\in\HH_\gengenus : \text{$\tr(x^2)\le c^{-1}$ and~$y\ge cI_\gengenus$}
	\}.
	\]
	\begin{prop}\label{prop;convpropKEseries}
	Let~$f\in S^k_{r,\lattice}$.
	If~$k>\gengenus+r+1$, then the Klingen Eisenstein series~$\klinser{k}{\gengenus}{r}{\lattice}(f)$ converges absolutely and uniformly on any vertical strip of positive height, and is a modular form of genus~$\gengenus$ and weight~$k$ with respect to~$\weil{\lattice}{\gengenus}$.
	\end{prop}

	\subsection{The Siegel operator}\label{sec;Siegelphi}
	Let~$0\le r\le\gengenus$.	
	\begin{defi}
	Let~$\disceltwo\in D_\lattice^{\gengenus-r}$.
	The \emph{Siegel~$\Phi^{\gengenus-r}_\disceltwo$ operator} is defined as
	\be\label{eq;vvSiegelphi}
	\Phi^{\gengenus-r}_\disceltwo f(\tau)
	=
	\sum_{\discel\in D_\lattice^{r}}
	\lim_{\lambda\to\infty}
	f_{(\discel,\disceltwo)}\left(\left(\begin{smallmatrix}
	\tau & 0\\
	0 & i\lambda I_{\gengenus-r}
	\end{smallmatrix}\right)\right)\mathfrak{e}_\discel,\qquad \tau\in\HH_{r},
	\ee
	for every~$f=\sum_{\discelthree\in D_\lattice^\gengenus}f_\discelthree \mathfrak{e}_\discelthree\in M^k_{\gengenus,\lattice}$.
	If~$\gengenus-r=1$, then we simply denote the operator as~$\Phi_\disceltwo$.
	\end{defi}
	As for the scalar-valued case, the limit defining the Siegel operator converges since the Fourier expansion of~$f$ converges uniformly on vertical strips of positive heights.
	
	Note that if~$q(\disceltwo)\neq 0$, then~$\Phi^{\gengenus-r}_\disceltwo$ maps every Siegel modular form to~$0$.
	Moreover, if~$\disceltwo=(\disceltwo_j)_j$ with~$\beta_j\in D_\lattice$, then~$\Phi_\beta^{\gengenus-r}$ may be written as the composition of Siegel operators~$\Phi^{\gengenus-r}_\beta=\Phi_{\beta_{\gengenus-r}}\circ\dots\circ\Phi_{\beta_1}$.
	\begin{lemma}\label{lemma;phivviswelldef}
	Let~$\disceltwo\in D_\lattice^{\gengenus-r}$.
	If~$f\in M^k_{\gengenus,\lattice}$, then~$\Phi^{\gengenus-r}_\disceltwo f\in M^k_{r,L}$.
	\end{lemma}
	\begin{proof}
	It is analogous to~\cite[Section~$5$, Proposition~$1$]{klingen}.
	\end{proof}
\begin{lemma}\label{lemma;FexpofimPhiop}
	Let~$f\in M^k_{\gengenus,\lattice}$.
	If~$\disceltwo\in D_\lattice^{\gengenus-r}$ is isotropic, then
	\[
	\Phi_\disceltwo^{\gengenus-r} f(\tau)=
	\sum_{\discel\in D_\lattice^{r}}
	\sum_{\substack{T\in q(\discel) + \halfint{r} \\ T\ge 0}} 
	c_{\big(\begin{smallmatrix}
	T & 0\\
	0 & 0
	\end{smallmatrix}\big)}(f_{(\discel,\disceltwo)}) \cdot e\big(
	\trace(T\tau)
	\big)\mathfrak{e}_\discel,\qquad\text{$\tau\in\textcolor{\newcor}{\HH_{\gengenus-r}}$.}
	\]
	In particular, a Siegel modular form in~$M^k_{\gengenus,\lattice}$ is a cusp form if and only if it lies in~$\ker\Phi_\disceltwo$ for every~$\disceltwo\in D_{\lattice}$.
	\end{lemma}
	\begin{proof}
	As in the scalar-valued case, the limit~$\lim_{\lambda\to\infty}
	f_{(\discel,\disceltwo)}\left(\left(\begin{smallmatrix}
	\tau & 0\\
	0 & i\lambda I_{\gengenus-r}
	\end{smallmatrix}\right)\right)$ extracts the part of the Fourier expansion of~$f_{(\discel,\disceltwo)}$ associated to singular matrices.
	This can be checked as follows.
	
	Without loss of generality we consider only the case of~$r=\gengenus-1$.
	Let~$\tau\in\HH_{\gengenus-1}$.
	The sequence of points~$\left(\begin{smallmatrix}
	\tau & 0\\
	0 & i\lambda
	\end{smallmatrix}\right)$ with~$\lambda$ large is contained in a vertical strip of positive height in~$\HH_{\gengenus}$, on which the Fourier series of~$f$ converges uniformly.
	We may then compute the limit term by term as
	\ba\label{eq;inproofeasyahuu}
	\Phi_\beta f(\tau)
	&=
	\sum_{\discel\in D_\lattice^{\gengenus-1}}
	\sum_{\substack{T\in q(\discel,\beta)+\halfint{\gengenus} \\ T\ge 0}}c_T(f_{(\discel,\beta)})
	e\big(\tr(T^*\tau)\big)
	\underbrace{\lim_{\lambda\to\infty} e\big(i\lambda T_*\big)}_{=\delta_{T_*,0}}\mathfrak{e}_\discel
	\ea
	where we denoted by~$T^*$ and~$T_*$ respectively the top-left~$(\gengenus-1)\times(\gengenus-1)$-block and the bottom-right entry of~$T$.
	It is easy to see that if~$T\ge0$ and~$T_*=0$, then~$T=\left(\begin{smallmatrix}
	T^* & 0\\
	0 & 0
	\end{smallmatrix}\right)$.
	This and~\eqref{eq;inproofeasyahuu} give the Fourier expansion of~$\Phi_\disceltwo f$.
	Clearly, it vanishes for all~$\beta$ if and only if~$f$ is a cusp form.
	\end{proof}
	\begin{thm}\label{thm;Phi0onkling}
	Let~$f\in S^k_{r,\lattice}$.
	If~$k>\gengenus+r+1$, then~$\Phi_0^{\gengenus-r} \klinser{k}{\gengenus}{r}{\lattice} (f) =f$.
Moreover, if~$\disceltwo\in D_\lattice^{\textcolor{\newcor}{\genus-r}}$ is non-zero, then~$\Phi_\disceltwo^{\textcolor{\newcor}{\genus-r}} \klinser{k}{\gengenus}{r}{\lattice} (f) =0$.
	\end{thm}
	\begin{proof}
	By Proposition~\ref{prop;convpropKEseries} the Klingen Eisenstein series converges uniformly on any vertical strip of positive height.
	We may then rewrite~$\big(\Phi_\disceltwo^{\gengenus-r} \klinser{k}{\gengenus}{r}{\lattice} (f)\big)(\tau)$ as
	\be\label{eq;underunifconv}
	\sum_{\discel\in D_L^{r}}\sum_{\gamma\in\klimet{\gengenus}{r}\backslash \Mp_{2\gengenus}(\ZZ)}
	\lim_{\lambda\to\infty}
	\phi\left(\big(\begin{smallmatrix}
	\tau & 0 \\
	0 & i\lambda I_{\gengenus-r}
	\end{smallmatrix}\big)\right)^{-2k}
	\Big(
	\weil{\lattice}{\gengenus}(\gamma)^{-1}
	\big(f\big(\textcolor{\newcor}{\big(\gamma\cdot\big(\begin{smallmatrix}
	\tau & 0 \\
	0 & i\lambda I_{\gengenus-r}
	\end{smallmatrix}\big)\big)^*}\big)
	\otimes \zeroe{\gengenus-r}\big)
	\Big)_{(\discel,\disceltwo)}
	\mathfrak{e}_\discel.
	\ee
	We extract from~\eqref{eq;underunifconv} the terms associated to~$\gamma\in\klimet{\gengenus}{r}$.
	Recall that if~$\gamma=(M,\phi)$ with~$M\in\klipar{\gengenus}{r}$, then~$\phi(\tilde\tau)
	={\det D_4} ^{1/2} \phi^*(\tilde\tau^*)$ for every~$\tilde\tau\in\HH_\gengenus$, where~$D_4$ is the lower-right~$(\gengenus-r)\times (\gengenus-r)$ block of~$M$.
	This implies that the sum of the terms associated to~$\gamma\in\klimet{\gengenus}{r}$ extracted from~\eqref{eq;underunifconv} is
	\begin{align*}
	\sum_{\discel\in D_L^{r}}
	&
	\lim_{\lambda\to\infty}
	\phi\left(\big(\begin{smallmatrix}
	\tau & 0 \\
	0 & i\lambda I_{\gengenus-r}
	\end{smallmatrix}\big)\right)^{-2k}
	\Big(
	\weil{\lattice}{\gengenus}(\gamma)^{-1}
	\big(f(\gamma^*\cdot\tau)
	\otimes \zeroe{\gengenus-r}\big)
	\Big)_{(\discel,\disceltwo)}
	\mathfrak{e}_\discel
	\\
	&=
	\sum_{\discel\in D_L^{r}}
	\det D_4^{-k} \phi^*(\tau)^{-2k}
	\Big(
	\weil{\lattice}{\gengenus}(\gamma)^{-1}
	\phi^*(\tau)^{2k}\weil{r}{\lattice}(\gamma^*)f(\tau)
	\otimes \zeroe{\gengenus-r}\Big)_{(\discel,\disceltwo)}
	\mathfrak{e}_\discel
	\\
	&=
	\sum_{\discel\in D_L^{r}}
	\Big(
	f(\tau)
	\otimes \zeroe{\gengenus-r}\Big)_{(\discel,\disceltwo)}
	\mathfrak{e}_\discel
	=
	\begin{cases}
	f(\tau), & \text{if~$\disceltwo=0$,}\\
	0, & \text{otherwise,}
	\end{cases}
	\end{align*}
	where to deduce the first to the last equality above we used that
	\bes
	\weil{\lattice}{\gengenus}(\gamma)
	(\mathfrak{e}_\discel \otimes\zeroe{\gengenus-r})
	=
	{\det D_4} ^{-k}\cdot
	\big(
	\weil{\lattice}{r}(\gamma^*)
	\mathfrak{e}_\discel
	\big)\otimes\zeroe{\gengenus-r}
	\qquad
	\text{for every~$\discel\in D_\lattice^{r}$, $\gamma\in\klimet{\gengenus}{r}$.}
	\ees
	As one can easily show following the idea of~\cite[p.~$68$, Proof of Proposition~$5$]{klingen}, the terms in~\eqref{eq;underunifconv} arising from~$\gamma\not\in\klimet{\gengenus}{r}$ are zero.
	\end{proof}
	
	\subsection{The Petersson inner product}
	We quickly review the Petersson inner product and its properties for vector-valued modular forms with respect to the Weil representation.
	\begin{defi}
	Let~$f_1,f_2\in M^k_{\genus,\lattice}$, at least one of which is a cusp form.
	The \emph{Petersson inner product} of~$f_1$ and~$f_2$ is
	\[
	\pet{f_1}{f_2}\coloneqq\int_{\Sp_{2\genus}(\ZZ)\backslash\HH_\genus}
	\langle f_1(\tau), f_2(\tau)\rangle
	\det y^k
	\frac{dx\, dy}{\det y^{\genus+1}}.
	\]
	\end{defi}
	The Petersson inner product \textcolor{\quickref}{converges} and makes~$S^k_{\genus,\lattice}$ a Hermitian space.
%
	As in the classical setting, Klingen Eisenstein series are orthogonal to cusp forms with respect to~$\petempty$.
	This may be proved as in~\cite[Proof of Proposition~$8$, Part~(ii)]{klingen}.
	\begin{prop}\label{prop;ortcusp&kling}
	Let~$0\leq r <\genus$ and let~$f\in S^k_{r,\lattice}$.
	Then~$\klinser{k}{\gengenus}{r}{\lattice}(f)$ is orthogonal to~$S^k_{\genus,\lattice}$ with respect to the Petersson inner product.
	\end{prop}
	\subsection{Poincaré series}\label{sec;poincser}
	In this section we provide vector-valued generalizations of Poincaré series associated to positive definite matrices.
	Further information regarding the scalar-valued case \textcolor{\newcor}{and the vector-valued genus~1 case} can be found in respectively~\cite[Section~$75$]{klingen} \textcolor{\newcor}{and~\cite[Section~1.2]{bruinierbook}}.

	We denote by~$\trans{\gengenus}$ the subgroup
	\be\label{eq;translinsympl}
	\trans{\gengenus}
	\coloneqq
	\left\{
	\left(\begin{smallmatrix}
	I_\gengenus & B \\
	0 & I_\gengenus
	\end{smallmatrix}\right)
	: B\in\Sym_{\gengenus}(\ZZ)
	\right\}
	\ee
	of translations in~$\Sp_{2\gengenus}(\ZZ)$ and \textcolor{\newcor}{by~$\mettrans{\gengenus}$ its metaplectic cover}.
	\begin{defi}
	Let~$\discel\in D_\lattice^\gengenus$ and let~$T\in q(\discel)+\halfint{\gengenus}$ be a \emph{positive definite} matrix.
	Assume that~$k>2\genus$.
	The weight~$k$ Poincaré series~$\poincv{\gengenus}{\discel}{T}{\lattice}$ is defined as
	\[
	\poincv{\gengenus}{\discel}{T}{\lattice}(\tau)
	=
	\frac{1}{2}\sum_{\gamma\in\mettrans{\gengenus}\backslash\Mp_{2\gengenus}(\ZZ)} e\big(\tr(T\gamma\cdot\tau)\big) \cdot \phi(\tau)^{-2k} \weil{\lattice}{\gengenus}(\gamma)^{-1}\mathfrak{e}_\discel,
	\qquad\tau\in\HH_\gengenus.
	\]
	\end{defi}
	The Weil representation~$\weil{\lattice}{\genus}$ is trivial on the principal subgroup of level~$N$ of~$\Sp_\genus(\ZZ)$, where~$N$ is the level of the lattice~$L$.
	Hence, the vector-valued \textcolor{\newcor}{Poincaré} series~$\poincv{\gengenus}{\discel}{T}{\lattice}$ is a finite combination of scalar valued \textcolor{\newcor}{Poincaré} series for that principal congruence subgroup.
	The absolute and local uniform convergence of~$\poincv{\gengenus}{\discel}{T}{\lattice}$ is therefore implied by the convergence of scalar-valued Poincaré series; see e.g.~\cite{christian}.
	With an analogous argument, we may deduce that~$\poincv{\gengenus}{\discel}{T}{\lattice}$ is a cusp form.
	We also remark that~$\poincv{\gengenus}{\discel}{T}{\lattice}=\poincv{\gengenus}{-\discel}{T}{\lattice}$.

    \textcolor{\quickref}{The following result generalizes a well-known property of scalar-valued Poincaré series, and can be proved similarly as in~\cite[Section~$6$]{klingen}.}
	\begin{lemma}\label{lemma;petSiegpoinandcusp}
	Let~$\discel\in D_\lattice^\gengenus$, let~$T\in q(\discel) + \halfint{\gengenus}$ be \emph{positive definite}, and let~$f\in S^k_{\gengenus,\lattice}$.
	If~$k>2\genus$, then
	\bas
	\bigpet{f}{\poincv{\gengenus}{\discel}{T}{\lattice}}= \constant_{k,\gengenus}\cdot \det T^{(\gengenus+1)/2 - k} \cdot c_T(f_\discel),
	\eas
	where~$\constant_{k,\gengenus}$ is the factor depending only on~$k$ and~$\gengenus$ given by
	\[
	\constant_{k,\gengenus}
	=
	\pi^{\gengenus(\gengenus-1)/4}(4\pi)^{\gengenus(\gengenus+1)/2 - \gengenus k}
	\prod_{\ell=1}^\gengenus \Gamma\Big(k - \frac{\gengenus + \ell}{2}\Big).
	\]
    \textcolor{\quickref}{In particular, the Poincaré series of weight~$k$ and genus~$\genus$ span the space of cusp forms~$S^k_{\genus,\lattice}$.}
	\end{lemma}
	\subsection{The space of almost cusp forms and its modular decomposition}\label{sec;mod_dec}
	Let~$k>2\gengenus$, and let~$\klinser{k}{\gengenus}{r}{\lattice}\colon S^k_{r,\lattice}\to M^k_{\gengenus,\lattice}$ be the linear operator lifting vector-valued Siegel cusp forms of genus~$r$ and weight~$k$ to the associated Klingen Eisenstein series of genus~$\gengenus$.
	We define
	\[
	\vklisub{k}{\gengenus}{r}{\lattice}
	\coloneqq
	\klinser{k}{\gengenus}{r}{\lattice}(S^k_{r,\lattice})
\qquad\text{for every~$0\leq r\leq \genus$.}
	\]
	If~$L$ is unimodular, then these are the subspaces of Klingen Eisenstein series introduced in~\cite[Section~$5$, Definition~$4$]{klingen}.
	We remark that~$\vklisub{k}{\gengenus}{0}{\lattice}=\CC \eiszero{k}{\genus}{L}$ and~$\vklisub{k}{\gengenus}{\genus}{\lattice}=S^k_{\genus,\lattice}$.
	
	By Theorem~\ref{thm;Phi0onkling} the restriction of the Siegel operator~$\Phi_0^{\genus-r}\colon \vklisub{k}{\gengenus}{r}{\lattice}\to S^k_{r,\lattice}$ is an isomorphism whose inverse is~$\klinser{k}{\gengenus}{r}{\lattice}$.
	The subspace~$\vklisub{k}{\gengenus}{r}{\lattice}$ is orthogonal to~$S^k_{\genus,\lattice}$ for every~$0\leq r < \genus$ by Proposition~\ref{prop;ortcusp&kling}.
	\begin{defi}
	The space of \emph{almost cusp forms}~$\almcusp{k}{\genus}{\lattice}$ is the subspace of~$M^k_{\genus,\lattice}$ defined as
	\[
	\almcusp{k}{\genus}{\lattice}
	=
	\bigoplus_{r=0}^\genus
	\vklisub{k}{\gengenus}{r}{\lattice}.
	\]
	The \emph{modular decomposition of~$f\in\almcusp{k}{\genus}{\lattice}$} is the rewriting of~$f$ with respect to the direct sum above.
	\end{defi}
	
	If~$g=1$, then the space of almost cusp forms is the same as
    in~\cite[Definition~$2.7$]{peterson} and~\cite[p.~$500$,~(2)]{bruiniermöller}.
	Since the subspace of Eisenstein series in~$M^k_{1,L}$ is in general of dimension greater than~$\dim \vklisub{k}{1}{0}{\lattice}=1$, the space~$\almcusp{k}{1}{\lattice}$ is usually smaller than~$M^k_{1,L}$.
	The same holds in higher genus.
	
	As we will see later, the cohomology class~$[\KMtheta{\genus}]$ of the Kudla--Millson theta function is an almost cusp form.
	It is the purpose of the present article to describe the modular decomposition of~$[\KMtheta{\genus}]$.
	This is given in Section~\ref{sec;gen2}.
	
	We define the map~$\vklin{k}{\genus}{r}{\lattice}\colon \almcusp{k}{r}{\lattice}\to \almcusp{k}{\genus}{\lattice}$ as the inverse of the restriction of~$\Phi^{\genus-r}_0$ to~$\vklisub{k}{\gengenus}{0}{\lattice}\oplus\cdots\oplus \vklisub{k}{\gengenus}{r}{\lattice}$.
	Concretely, it is defined on modular decompositions as
	\be\label{eq;defofKlinblamap}
	\vklin{k}{\genus}{r}{\lattice}\colon \sum_{j=0}^r
	\klinser{k}{r}{j}{\lattice}(f_j)
	\longmapsto
	\sum_{j=0}^{r}
	\klinser{k}{\gengenus}{j}{\lattice}(f_j),
	\qquad\text{$f_j\in S^k_{j,L}$.}
	\ee

	\subsection{Diagonal restrictions of Siegel Eisenstein series}\label{sec;diagonalrestr_vv}
	Böcherer proved in~\cite[Satz~$12$]{doubling} that certain diagonal restrictions of (essential subseries of) Siegel Eisenstein series can be written in terms of Klingen Eisenstein series arising from Poincaré series.
	The goal of this section is to extend Böcherer's result to the case of vector-valued Siegel modular forms with respect to the Weil representation.
	This is then used to show that the vanishing of Garrett's map stated in~\cite[Satz~$13$~$(1)$]{doubling} holds also in the vector-valued setting.
	This is a key ingredient to prove the modular decomposition of the Kudla--Millson theta function.
	\\
	
	
	Let~$r\leq\genus$ be a non-negative integer.
	We denote the \emph{Fourier--Jacobi expansion of cogenus~$r$ of~$\eiszero{k}{\gengenus}{\lattice}$} as
	\be\label{eq;FJexpsvcvv}
	\eiszero{k}{\gengenus}{\lattice}(\tau)
	=
	\sum_{\discel\in D_\lattice^r}
	\sum_{\substack{T\in q(\discel) + \halfint{r} \\ T\ge 0}}
	\mathfrak{e}_\discel
	\otimes
	\jaccoefvv{\genus}{r}{\discel}(\tau_4,\tau_2,T)\cdot e(\tr T\tau_1),
	\ee
	where~$\tau=\big(\begin{smallmatrix}
	\tau_1 & \tau_2\\
	\tau_2^t & \tau_4
	\end{smallmatrix}\big)$ is such that~$\tau_1\in\HH_r$ and~$\tau_4\in\HH_{\genus-r}$.
	The Fourier--Jacobi coefficients
	\[
	\jaccoefvv{\genus}{r}{\discel}\colon \HH_{\genus-r}\times\CC^{r\times(\genus-r)}\to\CC[D_\lattice^{\genus-r}]
	\]
	are Jacobi forms of weight~$k$ with respect to a Weil representation associated to the metaplectic Jacobi group of cogenus~$r$; see e.g.~\cite[Section~$3.2$]{kieferzuffetti} for the case of genus~$2$.
	\begin{defi}\label{def;Essentialpart}
	The \emph{essential part} of the Fourier--Jacobi expansion~\eqref{eq;FJexpsvcvv} of~$\eiszero{k}{\gengenus}{\lattice}$ is
	\[
	\espart{k}{\genus}{r}(\tau)
	=
	\sum_{\discel\in D_\lattice^r}
	\sum_{\substack{T\in q(\discel) + \halfint{r} \\ T>0}}
	\mathfrak{e}_\discel
	\otimes
	\jaccoefvv{\genus}{r}{\discel}(\tau_4,\tau_2,T)\cdot e(\tr T\tau_1).
	\]
	\end{defi}
	If~$M$ is a~$m\times n$ matrix with coefficients in~$\QQ$, in short~$M\in\QQ^{(m, n)}$, then we write the block decomposition of~$M$, where the top-left block is a~$\nu\times\nu$ matrix, as
	\[
	M=\Big(\begin{smallmatrix}
	M_1^{(\nu)} & M_2^{(\nu,n-\nu)}\\
	M_3^{(m-\nu,\nu)} & M_4^{(m-\nu,n-\nu)}
	\end{smallmatrix}\Big)\qquad\text{for every~$0\leq \nu \leq \min\{m,n\}$.}
	\]
	Let~$\widetilde{\mathfrak{A}}_{\genus,r}$ be the metaplectic cover of
	\begin{align*}
	\mathfrak{A}_{\genus,r}
	&=
	\left\{
	M=\big(
	\begin{smallmatrix}
	A & B \\
	C & D
	\end{smallmatrix}
	\big)\in\Sp_{2\genus}(\ZZ) : \text{$C=0$, $D=\big(\begin{smallmatrix}
	I_r & \ast \\
	0 & \ast
	\end{smallmatrix}\big)$, $A_1^{(r)}=I_r$}
	\right\}.
	\end{align*}
	\begin{defi}
	Let~$\discel\in D_\lattice^r$, and let~$t\in q(\discel) + \halfint{r}$ be a positive definite matrix.
	We define~$T\coloneqq\big(\begin{smallmatrix}
	t & 0\\
	0 & 0
	\end{smallmatrix}\big)\in q(\discel,0,\dots,0) + \halfint{\genus}$.
	For every~$k>\genus+r+1$, the \emph{generalized Poincaré series of exponential type}~$\mathfrak{g}_{\genus,r}^k(T,\discel)$ is
	\[
	\mathfrak{g}_{\genus,r}^k(\tau,T,\discel)
	\coloneqq
	\frac{1}{2}
	\sum_{\gamma\in\widetilde{\mathfrak{A}}_{\genus,r}\backslash\Mp_{2\genus}(\ZZ)}
	e(\tr T \gamma\cdot\tau)
	\phi(\tau)^{-2k}
	\weil{\lattice}{\genus}(\gamma)^{-1}
	\mathfrak{e}_\discel
	\otimes
	\zeroe{\genus-r},
	\qquad
	\tau\in\HH_\genus.
	\]
	\end{defi}

	If~$r=\genus$, then~$\mathfrak{g}_{\genus,\genus}^k(T,\discel)$ equals the Poincaré series~$\poincv{\gengenus}{\discel}{T}{\lattice}$ introduced in Section~\ref{sec;poincser}.
	We also remark that the generalized Poincaré series~$\mathfrak{g}_{\genus,r}^k(T,\discel)$ is actually the Klingen Eisenstein series arising from the Poincaré series~$\poincv{r}{\discel}{t}{\lattice}$, namely
	\[
	\mathfrak{g}_{\genus,r}^k(\tau,T,\discel)
	=
	\klinser{k}{\gengenus}{r}{\lattice}(\tau,\poincv{r}{\discel}{t}{\lattice}).
	\]
	Let
	\bas
	\ZZ^{(\genus-r,r)}_s
	&\coloneqq
	\{
	w_3\in\ZZ^{(\genus-r,r)}
	:
	\rk(w_3)=s
	\},
	\\
	\ZZ^{(\genus-r,r)}_{s,0}
	&\coloneqq
	\left\{
	w_3 \in\ZZ^{(\genus-r,r)}_s
	:
	w_3=\left(\begin{smallmatrix}
	(w_3)_1^{(s)} & \ast
	\\
	0 & 0_{\genus-r-s,r-s}
	\end{smallmatrix}\right)
	\right\},
	\\
	\GL(\genus,\ZZ)_\nu
	&\coloneqq
	\{
	U\in\GL(\genus,\ZZ) : U_3^{(\genus-\nu,\nu)}=0
	\},
	\\
	M^\genus_r(\ZZ)^*
	&\coloneqq
	\{
	W\in\ZZ^{(r,r)} : \text{$\det W\neq 0$ and there is~$V\in \GL(\genus,\ZZ)$ with~$V_1^{(r)}=W$}
	\}.
	\eas
	
	We are now ready to state the vector-valued generalization of Böcherer's result on diagonal restrictions of principal parts of Siegel Eisenstein series, which can be proved following the same wording of~\cite[Satz~$12$]{doubling}.
	We denote by~$A[B]$ the product of matrices~$B^t A B$.
	If~$A,B$ are integral matrices with the same number of columns, we say that~$\big(\begin{smallmatrix}
	A \\ B
	\end{smallmatrix}\big)$ is \emph{primitive} if that tuple can be extended to a square matrix~$\big(\begin{smallmatrix}
	A & \ast \\
	B & \ast
	\end{smallmatrix}\big)$ that is invertible over~$\ZZ$.
	\begin{thm}[Böcherer]\label{thm;genBöcherervv}
	If~$k>\genus+1$, $\tau_1\in\HH_r$ and~$\tau_4\in\HH_{\genus-r}$, then
	\bas
	\espart{k}{\genus}{r}
	\big(\begin{smallmatrix}
	\tau_1 & 0\\
	0 & \tau_4
	\end{smallmatrix}\big)
	&=
	\sum_{\discel\in D_\lattice^r}
	\sum_{\substack{T\in q(\discel) + \halfint{r} \\ T>0}}
	\sum_{w_1\in M^\genus_r(\ZZ)^*/\GL(r,\ZZ)}
	a_r^k(T,\discel)
	\cdot e\big(
	\tr T[w_1^t]\tau_1
	\big)
	\\
	&\quad\times
	\sum_{s=0}^{\min\{r,\genus-r\}}
	\varepsilon(s)
	\sum_{\substack{w_3\in\GL(\genus-r,\ZZ)_s\backslash\ZZ^{(\genus-r,r)}_{s,0} \\ \left(\begin{smallmatrix}
	w_1 \\ w_3
	\end{smallmatrix}\right)\text{ primitive}}}
	\mathfrak{e}_{\discel w_1^t}
	\otimes
	\mathfrak{g}_{\genus-r,s}^k\big(\tau_4,T[w_3^t],\discel(w_3^t)_{(s)}\big),
	\eas
	where~$(w_3^t)_{(s)}$ denotes the matrix given by the first~$s$ columns of~$w_3^t$ and
	\[
	\varepsilon(s)=\begin{cases}
	1, & \text{if $s=0$,}\\
	2, & \text{if $s>0$.}
	\end{cases}\]
	\end{thm}

	The proof of the following auxiliary lemma is left to the reader.
	\begin{lemma}\label{lemma;propFJcoefSEvv}
	Let~$k>\genus+1$,
	and let~$T_r\in q(\discel) + \halfint{r}$ be a positive semidefinite matrix, where~$\discel\in D_\lattice^r$.
	Assume that there exists~$U\in \GL(\genus,\ZZ)$ such that~$U^t T_r U = \big(\begin{smallmatrix}
	0 & 0\\
	0 & T_{r-1}
	\end{smallmatrix}\big)$ for some~$T_{r-1}\in q(\discel') + \halfint{r-1}$, where~$\discel U=(\discel'',\discel')$ with~$\discel'\in D_\lattice^{r-1}$ and~$\discel''$ is isotropic.
	Then
	\be\label{eq;symmetryFjaccoefdiffgenvv}
	\jaccoefvv{\genus}{r}{\discel}(\tau_4,0,T_{r})
	=
	\begin{cases}
	0 & \text{if~$\discel''\neq 0$,}
	\\
	\jaccoefvv{\genus-1}{r-1}{\discel'}(\tau_4,0,T_{r-1}) & \text{otherwise},
	\end{cases}
	\ee
	for all~$\tau_4\in\HH_{\genus-r}$.
	Moreover, we have~$\jaccoefvv{\genus}{r}{0}(\tau_4,0,0)
	=
	\eiszero{k}{\genus-r}{\lattice}(\tau_4)$.
	\end{lemma}
	The following result is the previously announced vector-valued generalization of the vanishing of Garrett's map.
	\begin{cor}\label{cor;vanishingGmapvv}
	Let~$k>\genus+1$ and let~$\genus=r+s$ be such that~$1\leq r<s$.
	Then
	\be\label{eq;vanishingfromBösatz12vv}
	\bigpet{f}{\eiszero{k}{\genus}{\lattice}\big(
	\begin{smallmatrix}
	\tau_1 &   \\
	  & \ast
	\end{smallmatrix}\big)}=0
	\qquad
	\text{for every~$f\in S^k_{s,\lattice}$ and~$\tau_1\in\HH_r$}.
	\ee
	Here we denote by~$\ast$ the variable of integration in~$\HH_s$.
	\end{cor}
	\begin{proof}
	
	Since~$\jaccoefvv{\genus}{r}{\discel}(\tau_4,0,T_{r})$ is a Siegel modular form of weight~$k$ and genus~$s$ with respect to~$\weil{\lattice}{s}$ for every~$T_r\in q(\discel) + \halfint{r}$, we deduce that
	\be\label{eq;frmBöcherer12svv}
	\bigpet{f}{\eiszero{k}{\genus}{\lattice}\big(\begin{smallmatrix}
	\tau_1 &  \\
	 & \ast
	\end{smallmatrix}\big)}
	=
	\sum_{\discel\in D_\lattice^r}
	\sum_{\substack{T_r\in q(\discel) + \halfint{r}\\ T_r\geq 0}}
	\bigpet{f}{\jaccoefvv{\genus}{r}{\discel}(\ast,0,T_{r})}\cdot e(-\tr T_r\overline{\tau_1})
	\mathfrak{e}_\discel.
	\ee
	By Theorem~\ref{thm;genBöcherervv} the diagonal restriction~$\espart{k}{\genus}{r}\big(\begin{smallmatrix}
	\tau_1 & \\
	 & \tau_4
	\end{smallmatrix}\big)$ is, with respect to~$\tau_4\in\HH_s$, a combination of Klingen Eisenstein series of genus~$s$ arising from cusp forms of genus at most~$r$.
	Since~$r<s$, we deduce from Proposition~\ref{prop;ortcusp&kling} that
	\[
	\bigpet{f}{\espart{k}{\genus}{r}\big(\begin{smallmatrix}
	\tau_1 & \\
	 & \ast
	\end{smallmatrix}\big)}=0\qquad\text{for every~$\tau_1\in\HH_r$ and~$f\in S^k_{s,\lattice}$.}
	\]
	This implies that
	\[
	\bigpet{f}{\eiszero{k}{\genus}{\lattice}\big(\begin{smallmatrix}
	\tau_1 & \\
	 & \ast
	\end{smallmatrix}\big)}
	=
	\bigpet{f}{\underbrace{\eiszero{k}{\genus}{\lattice} - \espart{k}{\genus}{r}\big(\begin{smallmatrix}
	\tau_1 & \\
	 & \ast
	\end{smallmatrix}\big)}_{\eqqcolon \tespart{k}{\genus}{r}\big(\begin{smallmatrix}
	\tau_1 & \\
	 & \ast
	\end{smallmatrix}\big)}},
	\]
	where
	\[
	\tespart{k}{\genus}{r}\big(\begin{smallmatrix}
	\tau_1 & \\
	 & \tau_4
	\end{smallmatrix}\big)
	=
	\sum_{\discel\in D_\lattice^r}
	\sum_{\substack{T_r\in q(\discel) + \halfint{r} \\
	T_r\geq 0 \\ T_r\not> 0}}
	\mathfrak{e}_\discel\otimes
	\jaccoefvv{\genus}{r}{\discel}(\tau_4,0,T_r) \cdot e(\tr T_r\tau_1).
	\]
	Hence we have
	\be\label{eq;frombös12vv}
	\bigpet{f}{\eiszero{k}{\genus}{\lattice}\big(\begin{smallmatrix}
	\tau_1 &  \\
	 & \ast
	\end{smallmatrix}\big)}
	=
	\sum_{\discel\in D_\lattice^r}
	\sum_{\substack{T_r\in q(\discel) + \halfint{r}\\ T_r\geq 0 \\ T_r \not>0}}
	\bigpet{f}{\jaccoefvv{\genus}{r}{\discel}(\ast,0,T_{r})}\cdot e(-\tr T_r\overline{\tau_1})
	\mathfrak{e}_\discel.
	\ee
	
	We now prove~\eqref{eq;vanishingfromBösatz12vv} by induction on~$r$ and~$\genus$.
	If~$r=1$, then~\eqref{eq;frombös12vv} simplifies to
	\[
	\bigpet{f}{\eiszero{k}{\genus}{\lattice}\big(\begin{smallmatrix}
	\tau_1 & \\
	 & \ast
	\end{smallmatrix}\big)}
	=
	\sum_{\discel\in D_\lattice}
	\bigpet{f}{\jaccoefvv{\genus}{1}{\discel}(\ast,0,0)}
	\mathfrak{e}_\discel,
	\]
	which vanishes by Lemma~\ref{lemma;propFJcoefSEvv}.
	
	Let now~$r>1$.
	Suppose that~\eqref{eq;vanishingfromBösatz12vv} is true for~$r$ and every genus~$\genus$ such that~$r < \genus - r$.
	By~\eqref{eq;frmBöcherer12svv} we deduce that $\bigpet{f}{\jaccoefvv{\genus}{r}{\discel}(\ast,0,T_r)}=0$ for every~$\discel\in D_\lattice^r$, every positive semidefinite~$T_r\in q(\discel) + \halfint{r}$, and every~$f\in S^k_{\genus-r,\lattice}$.
	We show that~\eqref{eq;vanishingfromBösatz12vv} is true also for~$r+1$ and every genus~$\genus$ such that~$r+1<\genus-r-1$.
	By~\eqref{eq;frombös12vv}, it is enough to show that if~$f\in S^k_{g-r-1,\lattice}$, then
	\be\label{eq;fsvvdsvv}
	\bigpet{f}{\jaccoefvv{\genus}{r+1}{\discel}(\ast,0,T_{r+1})}=0
	\ee
	for every~$\discel\in D_\lattice^{r+1}$ and every positive semidefinite~$T_{r+1}\in q(\discel) + \halfint{r+1}$ that is \emph{not positive definite}.
	With respect to the action of~$\GL(r+1,\ZZ)$ on~$\halfint{r+1}$, every such a matrix~$T_{r+1}$ is equivalent to~$\big(\begin{smallmatrix}
	0 & 0\\
	0 & T_r
	\end{smallmatrix}\big)$, for some positive semidefinite~$T_r\in\halfint{r}$.
	Let~$U\in \GL(r+1,\ZZ)$ be such that~$U^t T_{r+1} U = \big(\begin{smallmatrix}
	0 & 0\\
	0 & T_r
	\end{smallmatrix}\big)$, and write~$\discel U = (\discel'',\discel')$ for some isotropic~$\discel''\in D_\lattice$.
	Since~${\jaccoefvv{\genus}{r+1}{\discel}}(\tau_4,0,T_{r+1})$ equals either zero or~$\jaccoefvv{\genus-1}{r}{\discel'}(\tau_4,0,T_r)$ by Lemma~\ref{lemma;propFJcoefSEvv}, we may deduce~\eqref{eq;fsvvdsvv} by induction.
	\end{proof}
    
	\section{Orthogonal Shimura varieties and their cohomology}\label{sec;digressdiffgeo}
	Smooth orthogonal Shimura varieties are (not necessarily compact) complete Kähler manifolds.
	In this section we explain to what extent \textcolor{\tocheck}{the known properties of compact Kähler manifolds, see e.g.~\cite{huybrechts} for an overview,} are still satisfied by orthogonal Shimura varieties.
	We also introduce the Kudla--Millson theta function and special cycles, and describe the irreducible components of the latter in terms of primitive lattice vectors.

	\subsection{$L^2$-cohomology groups}\label{sec;isoofseveralcoho}
	Let~$(X,\omega)$ be a (not necessarily compact) complete Kähler manifold of dimension~$n$ and volume form~$\omega^n$.
	We introduce~$L^2$-forms and~$L^2$-cohomology following~\cite[Chapter~XIV, Section~$3$]{borelwallach}, \cite{dai} and~\cite{saperzucker}.

    Let~$A^r(X)$ be the space of complex-valued differential~$r$-forms on~$X$.
    We denote by~$H^r(X)$ and~$H^{p,q}(X)$ respectively the de Rham cohomology group of degree~$r$ with coefficients in~$\CC$ and the Dolbeault cohomology group of bidegree~$(p,q)$.
    If~$\alpha\in A^r(X)$ is closed, then we denote by~$[[\alpha]]$ its de Rham cohomology class.
    
    Let
\be\label{eq;definprodkdifforms}
	(\alpha,\beta)\coloneqq \int_X \alpha\wedge\ast\overline{\beta}\qquad\text{for every~$\alpha,\beta\in A^r(X)$,}
	\ee
	where~$\ast\colon A^r(X)\to A^{2n-r}(X)$ is the Hodge~$\ast$-operator associated to  \textcolor{\tocheck}{the Kähler metric of~$X$.}
	This is not an inner product, since it may diverge.
	The \emph{$L^2$-norm} of a form~$\eta\in A^r(X)$ is~$(\eta,\eta)$.
	We denote by~$\sqforms^r(X)$ the subspace of~$A^r(X)$ containing all differential forms~$\eta$ such that
	\bes
	(\eta,\eta)<\infty\qquad\text{and}\qquad(d\eta,d\eta)<\infty,
	\ees
	i.e.\ both~$\eta$ and its exterior derivative are \emph{square-integrable}.
	We refer to the elements of~$\sqforms^r(X)$ as~\emph{$L^2$-forms}.
	The~\emph{$L^2$-cohomology}~$\sqH^*(X)$ is the cohomology arising from the complex given by~$\sqforms^*(X)$.
	\textcolor{\tocheck}{Note that if~${\alpha,\beta\in \sqforms^r(X)}$, then~$(\alpha,\beta)$ converges by the Cauchy--Schwarz inequality.}
	
	\begin{rem}\label{rem;completehencestokes}
	Since~$X$ is \emph{complete}, the Stokes Theorem holds by~\cite{gaffney}, in particular the integral on~$X$ of a top-degree exact~$L^2$-form vanishes.
	This implies that~$(\cdot{,}\cdot)$ extends to an inner product on the~$L^2$-cohomology groups.
	\end{rem}
	
	The proof of the Hodge Theorem in the case of~$X$ compact is based on certain isomorphisms from cohomology groups to spaces of \emph{harmonic forms}.
	The latter forms can be defined also in the non-compact case and may be used to prove an~$L^2$-version of the Hodge decomposition, as we are going to recall.
	
	An \emph{$L^2$-harmonic form} on~$X$ is an~$L^2$-form annihilated by the usual Laplacian~$\Delta=dd^*+d^*d$ such that also its (adjoint) differential is square-integrable.
	We denote by~$\sqharm^r(X)$ the space of~$L^2$-harmonic forms of degree~$r$.
	For details regarding the following result we refer to~\cite[pp.~$5$-$6$]{dai} and~\cite[Sections~$2.4$ and~$2.5$]{noetherlefconj}.
	\begin{thm}[$L^2$-Hodge Theorem]
	Let~$X$ be a complete Kähler manifold.
	If the cohomology group~$\sqH^r(X)$ is finite dimensional, then it is isomorphic to~$\sqharm^r(X)$ and admits a Hodge decomposition
	\[
	\sqH^r(X)=\bigoplus_{p+q=r}\sqH^{p,q}(X).
	\]
	\end{thm}
	
	The inclusion~$\sqforms^r(X)\hookrightarrow A^r(X)$ induces a homomorphism
	\be\label{eq;homonotalwiso}
	\sqH^r(X)\to H^r(X)
	\ee
	which is in general neither injective nor surjective.
	As we will see later, if~$r$ is sufficiently small and~$X$ is an orthogonal Shimura variety, then~\eqref{eq;homonotalwiso} is an isomorphism.

%

	
	\subsection{Orthogonal Shimura varieties}
	Let~$L$ be an even indefinite lattice of signature~$(n,2)$, and let~$G=\OO(L\otimes\RR)$.
	Let~$\domain$ be one of the two connected components of the complex manifold
	\[
	\mathcal{D}=\{ z\in\lattice\otimes\CC : \text{$q(z)=0$ and~$(z,\bar{z})<0$}\}/\CC^*.
	\]
	Then~$\domain$ is (a model of) the Hermitian symmetric domain of type~IV arising from~$G$.
	
	Let~$\disk(\lattice)$ be the discriminant kernel of~$\OO(\lattice)$, namely the kernel of the natural homomorphism~$\OO(\lattice)\to\Aut(\lattice'/\lattice)$, and let~$\disk{}^+(\lattice)=\disk(\lattice)\cap \OO^+(\lattice)$, where~$\OO^+(\lattice)$ is the group of isometries of~$\lattice$ preserving the component~$\domain$.
	\begin{defi}
	An \emph{orthogonal Shimura variety} is a quotient of the form~$X=\Gamma\backslash\domain$ for some subgroup~$\Gamma\subseteq\disk{}^+(\lattice)$ of finite index.
	\end{defi}
	By the Theorem of Baily and Borel, the analytic space~$X$ admits a unique algebraic structure that makes it a quasi-projective variety.
	This may be compact only if~$n\le 2$.
	
	Assume that~$X$ is \emph{smooth}, e.g.\ by choosing~$\Gamma$ to be torsion-free.
	It is well-known that there is a unique~$G$-invariant Kähler form~$\omega$ on~$\domain$ up to a factor.
	Throughout this article, we choose~$\omega$ to be the Kähler form arising as the first Chern form of the (metrized) tautological line bundle of~$\domain$.
	The Hermitian metric on that bundle is the one induced by the bilinear form of~$L$; see~\cite[Section~$1$]{kudla-integrals} for details.
	Since~$X$ is smooth, the invariant form~$\omega$ induces a Kähler form on~$X$, which we still denote by~$\omega$.
	This makes~$X$ a complete Kähler manifold.
	
	
	Let~$\BBcomp{X}$ be the Baily--Borel compactification of~$X$, and let~$IH^*(\BBcomp{X})$ be its intersection cohomology.
	Zucker's conjectures, proved by Looijenga~\cite{looijenga} and Saper--Stern~\cite{saperstern}, asserts that there is an isomorphism~$\sqH^r(X)\cong IH^r(\BBcomp{X})$ for all~$r$.
	Moreover, since the dimension of the singular locus of~$\BBcomp{X}$ is at most~$1$, by~\cite[Theorem~$5.4$ and Remark~$5.5$~(ii)]{harriszucker} we also have an isomorphism~$IH^r(\BBcomp{X})\cong H^r(X)$ for all~$r<\dim(X)-1$.
	In fact, as illustrated in~\cite[Example~$3.4$]{noetherlefconj}, the map~\eqref{eq;homonotalwiso} is an isomorphism, and so is the composition
	\be\label{eq;3isoofcoho}
	\sqH^r(X)\cong IH^r(\BBcomp{X})\cong H^r(X)
	\qquad\text{for every~$r<\dim(X)-1$.}
	\ee
	From this we deduce that every de Rham cohomology class in~$H^r(X)$ admits a~$L^2$-representative.
	By~\cite[Theorem~$5.4$]{harriszucker} the composition of isomorphisms in~\eqref{eq;3isoofcoho} is a Hodge structure morphism.
	
	The Lefschetz operator
    \[
    \textcolor{\tocheck}{\lefschetz\colon H^{r}(X)\to H^{r+2}(X), \qquad
    [[\alpha]]\mapsto[[\alpha\wedge \omega]]}
    \]
    is injective for all $r<\dim(X)-2$, although~$X$ is only quasi-projective.
	This can be deduced in terms of the analogue of the Hard Lefschetz Theorem~\cite[Corollary~$9.2.3$]{maxim} for the intersection cohomology~$IH^r(\BBcomp{X})$.
	In fact, the Kähler class~$\omega$ of~$X$ is identified with the Chern class of an ample line bundle on~$\BBcomp{X}$; see~\cite[Sections~$2.4$ and~$2.5$]{noetherlefconj}.
	Moreover, the Lefschetz decomposition is compatible with the Hodge structure; see e.g.~\cite[Remark~$2.3$]{noetherlefconj}.
	
	We summarize the results above as follows.
    We denote the inner product on~$H^r(X)$ induced by~$(\cdot{,}\cdot)$ under~\eqref{eq;3isoofcoho} with the same symbol.
    \textcolor{\tocheck}{Let
     \bas
	 \sqHprim^r(X)
	 &\coloneqq
	 \ker\big(
	 \duallef\colon \sqH^r(X)\to \sqH^{r-2}(X)
	 \big),
	 \\
	 \sqHprim^{p,q}(X)
	 &\coloneqq
	 \ker\big(
	 \duallef\colon \sqH^{p,q}(X)\to \sqH^{p-1,q-1}(X)
	 \big)
	 \eas
    be the \emph{primitive~$L^2$-cohomology groups} of~$X$, where~$\Lambda\coloneqq\ast^{-1}\circ\mathcal{L}\circ\ast$ is the dual Lefschetz operator.}
    Similarly,~$\Hprim^r(X)$ and~$\Hprim^{p,q}(X)$ denote the primitive cohomology groups of~$X$.
	\begin{thm}\label{thm;lefdecL2}
	Let~$X$ be a smooth orthogonal Shimura variety of dimension~$n$.
	Let~$r<n-1$, so that the de Rham cohomology group~$H^{r}(X)$ is identified with~$\sqH^r(X)$ as in~\eqref{eq;3isoofcoho}.
	We have a decomposition
	\be\label{eq;lefdecL2inthm}
	H^r(X)=\bigoplus_{j\ge 0} \lefschetz^j \Hprim^{r-2j}(X),
	\ee
	which is orthogonal with respect to the inner product~$(\cdot{,}\cdot)$.
	Moreover, if~$r< n-2$, then the Lefschetz operator~$\lefschetz\colon H^r(X)\to H^{r+2}(X)$ is injective.
	Analogous results hold on each term of the Hodge decomposition of~$H^r(X)$.
	\end{thm}
	The fact that the Lefschetz decomposition of Theorem~\ref{thm;lefdecL2} is orthogonal with respect to~$(\cdot{,}\cdot)$ may be proved in the same way as for the compact case.
	In fact,~$L^2$-representatives of cohomology classes in different terms of the decomposition~\eqref{eq;lefdecL2inthm} are fiberwise orthogonal with respect to~$\langle\cdot{,}\cdot\rangle$ by e.g.~\cite[Proposition~$1.2.30$]{huybrechts}.
	
	\begin{rem}
	If~$\beta\in \lefschetz^j \sqHprim^{r-j,r-j}(X)$, then
	\be\label{eq;formmatchonHggsq}
	(\alpha,\beta)
	=
	\frac{(-1)^{(r-j)(2(r-j)+1)}\cdot n! \cdot j!}{(n-2r+j)!}
	\int_X\alpha\wedge\overline{\beta}\wedge\omega^{n-2r}
	\qquad
	\text{for every~$\alpha\in \sqH^{r,r}(X)$}.
	\ee
	This formula holds even if~$X$ is non-compact, and follows easily from e.g.~\cite[Proposition~$1.2.31$]{huybrechts}.
    The latter is a \emph{local} result and does not depend on the compactness of~$X$.
	The convergence of~$(\alpha,\beta)$ is ensured by the square-integrability.
	\end{rem}

	\subsection{Special cycles and the Kudla--Millson theta function}\label{sec;KMthetabasics}
	
	In the \textcolor{\newcor}{1980's} Kudla and Millson constructed a theta function~$\KMtheta{\gengenus}=\KMtheta{\gengenus}(\tau,z)$ in the two variables~$\tau\in\HH_\gengenus$ and~$z\in\domain$, with several remarkable properties.
	We provide here a quick overview of its construction together with the main results needed for the purposes of this article.
	We refer to~\cite{kudlamillson-harmonic1}, \cite{kudlamillson-harmonic2} and~\cite{kudlamillson-intersection} for further information.
	
	Let~$V_\RR=L\otimes\RR$.
	Kudla and Millson defined a~$G$-invariant Schwartz function~$\KMschwartz{\gengenus}$ on~$V_\RR^\gengenus$ with values in the space of smooth closed~$(\gengenus,\gengenus)$-forms on~$\domain$.
	This is first constructed on a base-point of~$\domain$ by applying a differential operator to the standard Gaussian of~$V_\RR^\gengenus$, and then spread to the whole~$\domain$ by pulling-back; see~\cite[Section~$3$]{kudlamillson-harmonic1} for further information.
	
	Let~$\schrod{\gengenus}$ be the Schrödinger model of (the restriction of) the Weil representation of~$\Mp_{2\gengenus}(\RR)$ acting on the space of Schwartz functions on~$V_\RR^\gengenus$; see e.g.~\cite[Section~$4$]{funkemillsonhyp} for a detailed definition.
	Let~$k=1+n/2$.
	The Kudla--Millson theta function of genus~$\gengenus$ associated to the lattice~$\lattice$ is defined as
	\[
	\KMtheta{\gengenus}(\tau,z)=\det y ^{-k/2}
	\sum_{\discel\in D_\lattice^\gengenus}
	\sum_{\lambda \in \discel + \lattice^\gengenus}
	\big(\schrod{\gengenus}(g_\tau)\KMschwartz{\gengenus}\big)(\lambda,z)
	\mathfrak{e}_\discel,
	\]
	for every~$\tau=x+iy\in\HH_\gengenus$ and~$z\in\domain$, where~$g_\tau$ is the standard element of~$\Sp_{2\gengenus}(\RR)$ mapping~$i\in\HH_\gengenus$ to~$\tau$.
	If~$\genus=0$, we define~$\KMtheta{0}$ to be the constant function equal to~$1$.
	
	With respect to the symplectic variable, the theta function~$\KMtheta{\gengenus}$ \textcolor{\newcor}{transforms} as a (non-holomorphic) Siegel modular form of weight~$k$ and genus~$\gengenus$ with respect to the Weil representation~$\weil{\lattice}{\gengenus}$.
	With respect to~$z\in\domain$, it is a~$\disk{}^+(\lattice)$-invariant closed~$(\gengenus,\gengenus)$-form on~$\domain$.
	Hence, the theta function descends to a function on the smooth orthogonal Shimura variety~$X=\Gamma\backslash\domain$, which we still denote by~$\KMtheta{\gengenus}$.
	
	An interesting feature of that theta function is that its cohomology class behaves as a \emph{holomorphic} Siegel modular form, since the image of~$\KMtheta{\gengenus}$ under~$\bar\partial$ in the symplectic variable~$\tau$ is an exact differential form on~$X$. 
	This means that the de Rham cohomology class~$\drcoho{\KMtheta{\gengenus}}$ may be regarded as an element of~$M^k_{\gengenus,\lattice}\otimes H^{\gengenus,\gengenus}(X)$.
	
	To describe the Fourier expansion of~$\drcoho{\KMtheta{\gengenus}}$ we need to introduce the so-called \emph{special cycles} of~$X$.
	For any vector~$\lambda\in (L')^\gengenus$ we denote by~$\lambda^\perp$ the set of points in~$\domain$ orthogonal to every component of~$\lambda$.
	The special cycle of~$X$ associated with~$\discel\in D_\lattice^\gengenus$ and~$T\in q(\discel)+\halfint{\gengenus}$, where~$T\ge 0$, is
	\be\label{eq;defspcy}
	Z(T,\discel) = \Gamma\big\backslash\bigg(
	\sum_{\substack{\lambda\in\discel + \lattice^\gengenus\\ q(\lambda)=T}}
	\lambda^\perp
	\bigg).
	\ee
	If non-empty, this is a codimension~$\rk (T)$ effective cycle on~$X$.
	We denote the cohomology class induced by~$Z(T,\discel)$ in~$H^{2\rk(T)}(X,\CC)$ by~$\drcoho{Z(T,\discel)}$.
	Note that~$Z(0,0)=X$ and that if~$\alpha\neq 0$, then~$Z(0,\alpha)=\emptyset$.
	Special cycles of codimension~$1$ are also known as \emph{Heegner divisors}.
	
	The group~$\GL_\gengenus(\ZZ)$ acts on~$\halfint{\gengenus}$ under the action~$T\mapsto M^t T M$, for~$T\in\halfint{\gengenus}$ and~$M\in\GL_\gengenus(\ZZ)$.
	This induces the symmetry
	\be\label{eq;symmetryspcy}
	Z(T,\discel) = Z(M^tTM,\discel M).
	\ee
	
	The cohomology class~$\drcoho{\KMtheta{\gengenus}}$ is the generating series of (the cohomology classes of) the codimension~$\gengenus$ special cycles, in fact the Fourier expansion of~$\drcoho{\KMtheta{\gengenus}}$ is
	\[
	\drcoho{\KMtheta{\genus}(\tau)}=\sum_{\discel\in D_\lattice^\gengenus}\sum_{\substack{T\in q(\discel) + \halfint{\gengenus}\\ T\ge0}}
	\drcoho{Z(T,\discel)}
	\wedge
	\drcoho{-\omega}^{\gengenus-\rk(T)}
	q^T \mathfrak{e}_\discel.
	\]
	
	\begin{rem}\label{rem;resandpullbKMthetafunct}
	Let~$M\subset\lattice$ be a sublattice of~$L$.
	Let~$\Gamma_L=\disk{}^+(\lattice)$ and let~$\Gamma_M=\disk{}^+(M)$.
	We denote by~$\KMtheta{\genus}^\lattice$, resp.~$\KMtheta{\genus}^M$, the Kudla--Millson theta function on~$X_L=\Gamma_L\backslash\domain$ arising from~$L$, resp.\ on~$X_M=\Gamma_M\backslash\domain$ arising from~$M$.
	It is easy to see that if~$\pi\colon X_M\to X_\lattice$ is the quotient map
	and if~$\tr_{\genus,\lattice/M}\colon M^k_{\genus,M}\to M^k_{\genus,\lattice}$ is as in Lemma~\ref{lemma;genofBYandS}, then
	\be\label{eq;resandpullbKMthetafunct}
	\pi^*\big(\drcoho{\KMtheta{\genus}^\lattice}\big)
	=
	\tr_{\genus,\lattice/M}\big(\drcoho{\KMtheta{\genus}^M}\big).
	\ee

	\end{rem}
	
	We conclude this section by recalling a product formula for the above theta functions; see e.g.~\cite[$(9.2)$]{kudla-algcycles}.
	As usual, we identify products of group algebras as in~\eqref{eq;identgroupalg}.
	Let~$r_1,r_2\in\ZZ_{\ge 0}$.
	If~$\tau_1\in\HH_{r_1}$ and~$\tau_2\in\HH_{r_2}$, then
	\be\label{eq;wedgeprodformforKMtheta}
	\KMtheta{r_1}(\tau_1)
	\wedge
	\KMtheta{r_2}(\tau_2)
	=
	\KMtheta{r_1+r_2}\big(
	\begin{smallmatrix}
	\tau_1 & 0\\
	0 & \tau_2
	 \end{smallmatrix}\big).
	\ee

	\subsection{Adelic and primitive special cycles}\label{sec;adelicspcy}
    \textcolor{\newcor}{In this section, we study the irreducible components of codimension~$g$ special cycles and address the question of whether they span the whole cohomology group~$H^{g,g}(X)$.
    This section may be skipped on a first reading.}
    
	Let~$V=\lattice\otimes\QQ$, and let~$\mathcal{S}\big(V^\genus(\Ab_f)\big)$ be the space of Schwartz functions on~$V^\genus(\Ab_f)$.
	For each prime~$p$, we let~$\lattice_p=\lattice\otimes\ZZ_p$ and let~$K_p$ be the subgroup of~$G(\QQ_p)$ which leaves~$\lattice_p$ stable and acts as the identity on~$\lattice_p'/\lattice_p$.
	Then~$K=\prod_p K_p$ is an open compact subgroup of~$G(\Ab_f)$.
	As above, we denote by~$\Gamma$ any finite index subgroup of~$K\cap G(\QQ)^+=\disk{}^+(L)$.
	
	\begin{defi}
	Given a~$K$-invariant Schwartz function~$\varphi\in\mathcal{S}(V^\genus(\Ab_f)\big){}^K$ and~$T\in\Sym_{\genus}(\QQ)$, we define the \emph{adelic special cycle}~$Z(T,\varphi,K)$ in~$X=\Gamma\backslash\domain$ as
	\[
	Z(T,\varphi,K)
	=
	\Gamma\big\backslash\sum_{\substack{x\in V^\genus \\ q(x)=T}}
	\varphi(x)\cdot x^\perp.
	\]
	\end{defi}
	
	Let~$\mathcal{S}_{L^\genus}\subset\mathcal{S}\big(V^\genus(\Ab_f)\big)$ be the subspace spanned by the characteristic functions~$\varphi_\lambda$ of the cosets~$\lambda + \widehat{\lattice^\genus}$ with~$\lambda\in (\lattice')^\genus$, where~$\widehat{\lattice^\genus}=\lattice^\genus\otimes\prod_p\ZZ_p$.
	The special cycles defined in Section~\ref{sec;KMthetabasics} are particular adelic special cycles, more precisely
	\[
	Z(T,\lambda + L^\genus)
	=
	Z(T,\varphi_\lambda,K)\qquad
	\text{for all~$\lambda\in (L')^\genus$.}
	\]
	
	We recall below the main result of~\cite{noetherlefconj} on the span of cohomology classes of adelic special cycles.
	
	\begin{thm}[Bergeron--Li--Millson--Moeglin]\label{thm;noetherlef}
	If~$\Gamma=\disk{}^+(\lattice)$ and~$\genus <(n+1)/3$, then~$H^{\genus,\genus}(X)$ is spanned by the classes of adelic special cycles of the form
	\[
	\bigdrcoho{Z(T,\varphi,K)}\wedge\bigdrcoho{\textcolor{\newcor}{(-\omega)}^{\genus-\rk T}},
	\qquad
	\text{$T\in\Sym_{\genus}(\QQ)$, $\varphi\in\mathcal{S}(V^\genus(\Ab_f)\big){}^K$.}
	\]
	\end{thm}
	It is natural to ask whether a result analogous to Theorem~\ref{thm;noetherlef} holds also for (non-adelic) special cycles.
	
	\begin{rem}\label{rem;BLMMhasimprecision}
	It is stated in~\cite[Remark~$2.4$]{noetherlefconj} that if~$\Gamma=\disk{}^+(L)$, then the adelic special divisors on~$X=\Gamma\backslash\domain$ are linear combinations of Heegner divisors.
	The argument provided in~\cite{noetherlefconj}, i.e.\ that every~$K$-invariant Schwartz function~$\varphi\in\mathcal{S}(V(\Ab_f))$ corresponds to a linear combination of characteristic functions on~$D_L$, is incorrect.
	Counterexamples are given by the characteristic function of~$N\widehat{L}$, with~$N\in\ZZ$, $N>1$.
	Hence, we can not deduce from~\cite{noetherlefconj} that~$H^{1,1}(X)$ is generated by (non-adelic) Heegner divisors.
	\end{rem}
	
	We illustrate below a non-adelic version of Theorem~\ref{thm;noetherlef}, proved in terms of primitive special cycles.
	For the sake of simplicity, we illustrate it for the cases of codimension~$1$ and~$2$.
	Together with the injectivity of the Kudla--Millson lift (known up to genus~$2$) and the main results of Section~\ref{sec;gen2}, the non-adelic version of Theorem~\ref{thm;noetherlef} will imply dimension formulas for cohomology groups of orthogonal Shimura varieties; see Corollaries~\ref{cor;dimcoho1} and~\ref{cor;dimcoho2}.

    \textcolor{\correction}{Recall from e.g.~\cite[Definition 2.17]{kneser} that a sublattice~$M\subseteq L'$ is said to be \emph{primitive in~$L'$} if there exists a basis of~$M$ which can be extended to a basis of~$L'$.
    By~\cite[Satz~14.3]{kneser}, $M$ is primitive in~$L'$
    if and only if~$M=(M\otimes\QQ)\cap L'$.}
    
    \begin{defi}
	\textcolor{\correction}{Let~$\alpha\in D_\lattice^\genus$ and let~$\primvec{\lattice}{\genus}{\alpha}$ be the set of tuples~$\lambda\in\alpha + \lattice^\genus$ such that the entries of~$\lambda$ span a primitive sublattice of~$L'$.
	For all~$T\in q(\alpha)+\halfint{\genus}$, the \emph{primitive special cycle}~$\Zprim(T,\alpha)$ is defined as}
	\[
	\textcolor{\newcor}{\Zprim(T,\alpha)
	=
	\Gamma\big\backslash\bigg(
	\sum_{\substack{\lambda\in\primvec{\lattice}{\genus}{\alpha}\\ q(\lambda)=T}}
	\lambda^\perp
	\bigg).}
	\]
    \end{defi}
    \begin{rem}\label{rem;ZpriminvGL}
    \textcolor{\correction}{If~$M\in\GL_\genus(\ZZ)$, then~$\Zprim(M^t T M,\alpha M)=\Zprim(T,\alpha)$.
    In fact, since the lattice spanned by the entries of~$\lambda$ equals the one spanned by the entries of~$\lambda M$, then the set~$\primvec{\lattice}{\genus}{\alpha}$ maps to~$\primvec{\lattice}{\genus}{\alpha M}$ under the right-multiplication by~$M$.
    Since~$M^t T M= q(\lambda M)$, we deduce that}
    \[
    \textcolor{\newcor}{\Zprim(T,\alpha)
    =
    \Gamma \big\backslash\bigg(
	\sum_{\substack{\lambda\in\primvec{\lattice}{\genus}{\alpha}\\ q(\lambda)=T}}
	(\lambda M)^\perp
	\bigg)
    =
    \Gamma \big\backslash\bigg(
	\sum_{\substack{\lambda\in\primvec{\lattice}{\genus}{\alpha} M\\ q(\lambda)=M^t T M}}
	\lambda^\perp
	\bigg)
    =
    \Zprim(M^t T M,\alpha M).}
    \]
    \end{rem}

    \textcolor{\correction}{We denote by~$M_\genus(\ZZ)^*$ the set of~$\genus\times\genus$ integral matrices with non-zero determinant.}
    \textcolor{\correction}{We say that a matrix~$A\in M_\genus(\ZZ)^*$ divides~$T\in q(\alpha)+\halfint{\genus}$, in short~$A|T$, if there exists~$\beta\in D_\lattice^\genus$ such that~$A^{-t} T A^{-1}$ lies in~$q(\beta) + \halfint{\genus}$.}

    \begin{defi}
    \textcolor{\newcor}{Let~$\mu\colon \GL_g(\ZZ)\backslash M_g(\ZZ)^*\to\ZZ$ be the \emph{generalized Möbius function} defined as~$\mu([A])= 1$ if~$A\in\GL_g(\ZZ)$, and as
    \[
    \mu([A])=
        -\sum_{\substack{[B]\in\GL_g(\ZZ)\backslash M_\genus(\ZZ)^* \\ \,B\not\in\GL_\genus(\ZZ) \\ AB^{-1}\in M_\genus(\ZZ)^*}}\mu([AB^{-1}])
    \]
    otherwise.}
    \end{defi}
    \textcolor{\newcor}{It is easy to see that the sum defining~$\mu$ is finite.
    Note that if~$g=1$, then $\GL_1(\ZZ)\backslash M_1(\ZZ)^*$ is isomorphic to~$\ZZ_{>0}$ and~$\mu$ boils down to the classical Möbius function.
    We also remark that}
    \begin{equation}\label{eq;muTU}
        \textcolor{\newcor}{\mu([A])=\mu([AU])\qquad \text{for all $U\in\GL_g(\ZZ)$.}}
    \end{equation}

    \begin{lemma}\label{lemma;spcyclesandprim}
        \textcolor{\correction}{Let~$\alpha\in D_\lattice^\genus$ and let~$T\in q(\alpha) + \halfint{\genus}$ be positive definite.
        Then}
	\be\label{eq;spcyclesandprim}
	\textcolor{\newcor}{Z(T,\alpha)
	=
	\sum_{\substack{[A]\in \GL_\genus(\ZZ)\backslash M_\genus(\ZZ)^*
	\\
	A|T}}
    \sum_{\substack{\beta\in D_\lattice^\genus \\ \beta A=\alpha}}
	\Zprim(A^{-t} T A^{-1},\beta)}
	\ee
	\textcolor{\newcor}{and}
	\be\label{eq;spcyclesandprim2}
	\textcolor{\newcor}{\Zprim(T,\alpha)
	=
	\textcolor{\newcor}{\sum_{\substack{[A]\in\GL_\genus(\ZZ)\backslash M_\genus(\ZZ)^* \\ A|T}}
    \mu([A])
    \sum_{\substack{\beta\in D_L^\genus \\ \beta A=\alpha}}
    Z(A^{-t}TA^{-1},\beta).}}
	\ee
    \end{lemma}
    \textcolor{\newcor}{We could drop the condition~$A|T$ in~\eqref{eq;spcyclesandprim} and~\eqref{eq;spcyclesandprim2}, since the cycles of index~$A^{-t}TA^{-1}$ are empty otherwise.}
    
    \begin{proof}
        \textcolor{\correction}{Let~$\lambda\in \alpha + L^\genus$ be with~$q(\lambda)=T$, and let~$S$ be the rank~$g$ sublattice of~$L'$ spanned by the entries of~$\lambda$.
        Then~$\lambda^\perp\subset\domain$ can be written as~$(S^{\operatorname{prim}}\otimes\QQ)^\perp$ for a unique primitive sublattice~$S^{\operatorname{prim}}$ of~$L'$.
        Let~$A_\lambda\in \GL_\genus(\ZZ)\backslash M_\genus(\ZZ)^*$ be the unique matrix dividing~$T$ such that the entries of~$\lambda A_\lambda^{-1}$ span~$S^{\operatorname{prim}}$.
        We may rewrite the formal sum of cycles in~$\domain$ defining~$Z(T,\alpha)$ as
        \[
       \sum_{\substack{\lambda\in \alpha + L^\genus \\ q(\lambda)=T}}\lambda^\perp
       =
       \sum_{\substack{\lambda\in \alpha + L^\genus \\ q(\lambda)=T}}(\lambda A_\lambda^{-1})^\perp
        =
        \sum_{\substack{A\in \GL_\genus(\ZZ)\backslash M_\genus(\ZZ)^*
	\\
	A|T}}
        \sum_{\substack{\beta\in D_L^\genus \\ \beta A=\alpha}}
    \sum_{\substack{\lambda\in\primvec{\lattice}{\genus}{\beta}\\ q(\lambda)=A^{-t}TA^{-1}}}
        \lambda^\perp,
        \]
        from which we deduce~\eqref{eq;spcyclesandprim}.
The decomposition~\eqref{eq;spcyclesandprim2} follows by an easy inductive argument on the number of classes~$[A]\in\GL_\genus(\ZZ)\backslash M_\genus(\ZZ)^*$ dividing~$T$.}\qedhere
    \end{proof}

	We denote by~$E_8$ the \textcolor{\newcor}{root lattice of the~$E_8$ root system}, and by~$U$ the hyperbolic plane.
	They are even unimodular lattices of signature~$(8,0)$ and~$(1,1)$ respectively.
	\begin{lemma}\label{lemma;orbitsofdisck}
	Let~$L=U^{\oplus 2}\oplus E_8\oplus E$ for some positive definite even lattice~$E$.
	There exists a constant~$C=C(E)\in\QQ_{>0}$ such that the following property is satisfied for all~$\alpha\in D_L^2$.
	If~$a,b\in \textcolor{\correction}{\primvec{\lattice}{2}{\alpha}}$ are such that~$q(a)=q(b)>0$ and~$q(a_1)>C$, then~$a$ and~$b$ lie in the same~$\disk{}^+(L)$-orbit.
	\end{lemma}
	\begin{proof}
	We adapt the argument of~\cite[Proof of Lemma~$4.4$]{freitaghermann}.
	Let
	\[
	L=U_1 \oplus U_2\oplus E_8\oplus E,
	\]
	where~$U_1$ and~$U_2$ are copies of the hyperbolic plane.
	
	There exists an isomorphism~$U_1\oplus U_2\to \Mat_2(\ZZ)$ such that the quadratic form of~$U_1\oplus U_2$ corresponds to the determinant on~$\Mat_2(\ZZ)$.
	The action of~$\SL_2(\ZZ)$ on~$\Mat_2(\ZZ)$ by multiplication from both sides induces isometries in~$\OO^+(U_1\oplus U_2)$.
	The theorem of elementary divisors for~$\SL_2(\ZZ)$ implies that every element of~$L'$ maps to~$L_1'$ under an isometry in~${\OO^+(U_1\oplus U_2)}$, where~$L_1\coloneqq \{0\}\oplus U_2\oplus E_8\oplus E$.
	
	Let
	\[
	C_\alpha\coloneqq \min\{
	q(\lambda) : \text{$\lambda\in \alpha + E$}
	\},
	\qquad
	\alpha\in D_{E}\cong D_{L}.
	\]
	Let~$C=\max_{\alpha} C_\alpha$, and let~$a,b\in(\lattice')^2$ be as in the statement of Lemma~\ref{lemma;orbitsofdisck}.
	Since~$q(a_1)>C$, there exists~$e\in\alpha+E$ such that~$q(e)<q(a_1)$.
	Since~$E_8$ represents primitively every positive integer, there exists a primitive~$\tilde e\in E_8$ such that~$q(\tilde e)=q(a_1)-q(e)$.
	This implies that~$\tilde e + e$ is primitive in~$E_8\oplus E'$ (hence in~$L'$), and has the same norm and same~$D_L$-coset as~$a_1$.
	By~\cite[Lemma~$4.4$]{freitaghermann} there exists an isometry in~$\disk{}^+(L)$ mapping~$a_1$ to~$\tilde e + e$.
	The same argument works for~$b_1$ in place of~$a_1$.
	We may therefore assume that~$a_1=b_1\in E_8\oplus E'$.

	Since~$E_8\oplus E'$ is preserved under the action of~$\OO^+(U_1\oplus U_2)$, by the theorem of elementary divisors for~$\SL_2(\ZZ)$ we may further assume that~$a_2,b_2\in L_1'$.
	To conclude the proof it is then enough to show that there exists~$\gamma\in\Stab_\Gamma(a_1)$ such that~$\gamma(a_2)=b_2$.
	We construct~$\gamma$ in terms of Eichler transformations.
	Recall that if~$u,v\in L$ are such that~$v\perp u$ and~$u$ is isotropic, then the Eichler transformation~$E(u,v)$ defined as
	\[
	E(u,v)(x)= x - (x,u) v + (x, v) u - q(v) (x,u) u,\qquad x\in L\otimes\QQ,
	\]
	lies in the discriminant kernel of~$L$.
	
	Let~$f_1,f_2$, resp.~$f_3,f_4$, be the standard generators of~$U_1$, resp.~$U_2$.
	Since~$a_1$ and~$a_2$ \textcolor{\correction}{span a primitive sublattice of~$L_1'$}, there exists a basis~$\mathcal{B}$ of~$L_1'$ whose first two vectors are~$a_1$ and~$a_2$.
	Let~$a'\in L_1$ be the negative of the second vector of the basis of~$L_1$ dual to~$\mathcal{B}$.
	Then
	\[
	(a_1,a')=0
	\qquad
	\text{and}
	\qquad
	(a_2,a')=-1
	.
	\]
	Similarly, one can find~$b'\in L_1$ such that~$(b_1,b')=0$ and~$(b_2,b')=-1$.
	It is easy to check that
	\[
	E(f_1,a')(a_i)=
	\begin{cases}
	a_i, & \text{if~$i=1$,}\\
	a_i-f_1, & \text{if~$i=2$,}
	\end{cases}
	\qquad\text{and}
	\qquad
	E(f_1,b')(b_i)=
	\begin{cases}
	b_i, & \text{if~$i=1$,}\\
	b_i-f_1, & \text{if~$i=2$.}
	\end{cases}
	\]
	We may further compute
	\bas
	E(f_2,b_2-a_2)(a_2-f_1)&=b_2-f_1,
	\\
	E(f_2,b_2-a_2)(a_1)&=a_1 + (a_1,b_2-a_2)f_2 = a_1,
	\eas
	where in the last equality we used that since~$q(a)=q(b)$, then
	\[
	(a_1,b_2-a_2)
	=
	(a_1,b_2)-(a_1,a_2)
	=
	(b_1,b_2)-(a_1,a_2)
	=0.
	\qedhere
	\]
	\end{proof}
	By~\cite[Lemma~$4.3$]{bruiniermöller}, if~$L$ splits off two hyperbolic planes, then the induced primitive special divisors on~$\disk{}^+(L)\backslash\domain$ are irreducible.
	The following result is a generalization of this in codimension~$2$.
	\begin{prop}\label{prop;onirredofprimcy}
	If~$X=\disk{}^+(\lattice)\backslash\domain$ and~$L$ splits off~$U^{\oplus 2}\oplus E_8$, then the primitive special cycles of codimension~$2$ arising from~$L$ are (not necessarily reduced) irreducible cycles.
	\end{prop}
	\begin{proof}
	Let~$\discel\in D_L^{2}$, and let~$T\in q(\alpha)+\halfint{2}$ be positive definite.
	We may choose~$M\in\GL_2(\ZZ)$ such that the top-left entry of~$M^t T M$ is greater than~$C$, where~$C$ is as in Lemma~\ref{lemma;orbitsofdisck}.
	By the latter lemma, the irreducible components of
	\[
	\sum_{\substack{\lambda\in\primvec{\lattice}{2}{\discel M}\\ q(\lambda)=M^tTM}}
	\lambda^\perp
	\]
	on~$\domain$ are identified under~$\disk{}^+(L)$.
	This means that the induced cycle~$\Zprim(M^tTM,\discel M)$ on~$X$ is irreducible.
	Since~$\Zprim(M^tTM,\discel M)=\Zprim(T,\alpha)$ \textcolor{\correction}{by Remark~\ref{rem;ZpriminvGL}}, also $\Zprim(T,\alpha)$ is irreducible.
	\end{proof}

	\begin{cor}\label{cor;blmminnonadcase}
	Let~$X=\disk{}^+(\lattice)\backslash\domain$ be of dimension~$n>2$.
	\begin{enumerate}
	\item If~$L$ splits off~$U^{\oplus 2}$, then~$H^{1,1}(X)$ is spanned by the classes of special divisors
	\[
	\bigdrcoho{Z(m,\alpha)}\wedge\bigdrcoho{\omega^{1-\delta_{0,m}}},
	\qquad
	\text{$\alpha\in D_\lattice$, $m\in q(\alpha) + \ZZ$, $m\ge0$.}
	\]
	\item If~$L$ splits off~$U^{\oplus 2}\oplus E_8$, then~$H^{2,2}(X)$ is spanned by the classes of special cycles\label{item;spcy2cod}
	\bes
	\bigdrcoho{Z(T,\alpha)}\wedge\bigdrcoho{\omega^{2-\rk T}},
	\qquad
	\text{$\alpha\in D_\lattice^2$, $T\in q(\alpha) + \halfint{2}$, $T\ge0$.}
	\ees
	\end{enumerate}
	\end{cor}
	\begin{proof}
	We give the details for the case of codimension~$2$, since a similar proof works also in codimension~$1$.
	By Proposition~\ref{prop;onirredofprimcy} and~\cite[Lemma~$4.3$]{bruiniermöller} the primitive special cycles of codimension at most~$2$ are irreducible.
	By Lemma~\ref{lemma;spcyclesandprim} the cohomology classes
	\[
	\bigdrcoho{\Zprim(T,\alpha)}\wedge\bigdrcoho{\omega^{2-\rk T}}
	\qquad
	\text{with $\alpha\in D_\lattice^2$ and $T\in q(\alpha) + \halfint{2}$, $T\ge 0$,}
	\]
	span the same subspace of~$H^{4}(X)$ spanned by the classes of special cycles with the same indices.
	To prove Corollary~\ref{cor;blmminnonadcase}~\eqref{item;spcy2cod}, it is enough to prove that the irreducible components of the adelic special cycles of codimension~$2$ are primitive special cycles.
	
	Let~$\varphi\in\mathcal{S}\big(V^2(\Ab_f)\big){}^K$.
	Then there exists an even lattice~$M\subset V$ such that~$\varphi\in\mathcal{S}_{M^2}$; see e.g.~\cite[\textcolor{\newcor}{fifth} line of the proof of Theorem~$3.2$]{kudla-speccyc}.
	The irreducible components of the adelic special cycle~$Z(T,\varphi,K)$ are then of the form
	\be\label{eq;adspircomp}
	\Gamma\big\backslash\sum_{x\in\Gamma \tilde{x}} x^\perp
	\ee
	for some~$\tilde{x}\in V^2(\QQ)\cap\supp(\varphi)$ with~$q(\tilde{x})=T$.
	Let~$m\in\ZZ_{>0}$ be such that~$\tilde\lambda\coloneqq m\tilde{x}\in(\lattice')^2$.
	Since~$\Gamma$ acts trivially on~$D_L^2$, there is some~$\alpha\in D_\lattice^2$ such that~$\Gamma\tilde\lambda\subseteq\alpha + \lattice^2$.
	Without loss of generality, we may assume that \textcolor{\correction}{the entries of~$\tilde\lambda$ span a primitive sublattice of}~$L'$, hence the subvariety~\eqref{eq;adspircomp} is (up to multiplicity) a primitive special cycle of coset~$\alpha$.
%
%
	\end{proof}
    
	\section{The generating series of volumes of special cycles}\label{sec;volofspeccyinKM}
	In this section we prove that~$\KMtheta{\genus}$ is square-integrable, and recall some integral formulas to compute volumes of special cycles.
	The latter are well-known for compact orthogonal Shimura varieties~\cite{kudlamillson-harmonic1} \cite{kudlamillson-harmonic2}, and have been extended to the non-compact case in~\cite{kudlamillson-tubes} and~\cite{kudlamillson-intersection}.
	\\
	
	We denote by~$X$ a smooth orthogonal Shimura variety arising from an indefinite even lattice~$L$ of signature~$(n,2)$.
	Let~$\discel\in D_\lattice^\gengenus$ and let~$T\in q(\discel)+\halfint{\gengenus}$ be positive definite.
	The differential form
	\[
	\zeta(T,\discel) \coloneqq 
	\exp\big(2\pi\trace T\big)
	\sum_{\substack{\lambda\in \discel+L^\gengenus \\ q(\lambda)=T}}
	\KMschwartz{\gengenus}(\lambda)
	\]
	is a Poincaré dual form of the special cycle~$Z(T,\discel)$; see~\cite[Section~$5$]{kudlamillson-tubes}.
	For any quadratic space~$V$, we denote by~$\Witt(V)$ its Witt index.
	\begin{lemma}\label{lemma;L2cohoKMtheta&fexp}
	If~$2\gengenus < n + 1 - \Witt(L\otimes\QQ)$, then the Poincaré dual form~$\zeta(T,\discel)$ of the special cycle~$Z(T,\discel)$ and the Kudla--Millson theta function~$\KMtheta{\gengenus}$ are~$L^2$-forms.
	If moreover~$2\genus < n-1$, then the~$L^2$-cohomology class~$[\KMtheta{\gengenus}]$ is a \emph{holomorphic} Siegel modular form of genus~$\gengenus$ and weight~$k=1+n/2$ with respect to~$\weil{L}{\genus}$ with Fourier expansion
	\be\label{eq;FexpKMthetaL2}
	[\KMtheta{\gengenus}]
	=
	\sum_{\discel\in D_\lattice^\gengenus}\sum_{\substack{T\in q(\discel) + \halfint{\gengenus}\\ T\ge 0}}
	[\zeta(T,\discel)]\wedge[-\omega]^{\gengenus-\rk(T)} q^T \mathfrak{e}_\discel.
	\ee
	\end{lemma}
	\begin{proof}
	To prove that~$\KMtheta{\gengenus}$ is square-integrable, we follow the argument of~\cite[Proof of Proposition~$4.1$]{bruinierfunke}. 
	It is enough to show that the~$\discel$-component of~$\KMtheta{\gengenus}$ is square-integrable for every~$\discel\in D_\lattice^\gengenus$.
	
	Let~$V=\lattice\otimes\RR$.
	The~$\discel$-component of~$\KMtheta{\gengenus}$ equals the theta function defined as
	\[
	\theta_\discel(\tau,\KMschwartz{\gengenus})
	\coloneqq
	\det y^{-k/2}\sum_{\lambda\in\discel + \lattice^{\gengenus}} \omega_\infty(g_\tau)\KMschwartz{\gengenus}(\lambda),
	\]
	where~$\omega_\infty$ denotes the Schrödinger model of (the restriction of) the Weil representation of~$\Mp_{2\gengenus}(\RR)$
acting on the Schwartz functions on~$V^\gengenus$.
	The square of the~$L^2$-norm of that differential form is
	\[
	\lVert\theta_\discel(\tau,\KMschwartz{\gengenus})\rVert^2_2
	=
	\int_X \theta_\discel(\tau,\KMschwartz{\gengenus}) \wedge \ast\overline{\theta_\discel(\tau,\KMschwartz{\gengenus})},
	\]
	where~$\ast$ is the Hodge star operator.
	Since~$\KMschwartz{\gengenus}$ is a \emph{real} differential form,
	it is easy to check that~$\overline{\omega_\infty(g_\tau)\KMschwartz{\gengenus}}=\omega_\infty(g_{-\overline{\tau}})\KMschwartz{\gengenus}$.
	We may then compute
	\[
	\ast\overline{\theta_\discel(\tau,\KMschwartz{\gengenus})}
	=
	\overline{\theta_\discel(\tau,\ast\KMschwartz{\gengenus})}
	=
	\theta_\discel(-\overline{\tau},\ast\KMschwartz{\gengenus}).
	\]
	Let~$\underline{\discel}=(\discel,\discel)\in D_\lattice^{2\gengenus}$ and let~$\KMschwartznew{\gengenus}(\underline{v})\nu=\KMschwartz{\gengenus}(v_1)\wedge\ast\KMschwartz{\gengenus}(v_2)$ for every~$\underline{v}=(v_1,v_2)\in V^2$, where~$\nu=\omega^n$ is the volume form of~$X$.
	Then
	\[
	\theta_\discel(\tau,\KMschwartz{\gengenus}) \wedge \ast\overline{\theta_\discel(\tau,\KMschwartz{\gengenus})}
	=
	\theta_{\underline{\discel}}(\tau,-\overline{\tau},\KMschwartznew{\gengenus})\nu
	\]
	may be regarded as the restriction to the diagonal of a theta function of genus~$2\gengenus$, where~$(\tau,-\overline{\tau})$ is identified with the element~$\left(\begin{smallmatrix}
	\tau & 0\\
	0 & -\overline{\tau}
	\end{smallmatrix}\right)\in\HH_{2\gengenus}$, and hence
	\be\label{eq;doublingsqint}
	\lVert\theta_\discel(\tau,\KMschwartz{\gengenus})\rVert^2_2
	=
	\int_X \theta_{\underline{\discel}}(\tau,-\overline{\tau},\KMschwartznew{\gengenus})\nu.
	\ee
	This integral converges by the Weil convergence criterion~\cite{weil}, since~$2\gengenus < n + 1 - \Witt(V)$.
	
	The proof that~$\zeta(T,\discel)$ is square-integrable is similar.
	We remark that
	\[
	\zeta(T,\discel)\wedge\ast\zeta(T,\discel)
	=
	\exp(4\pi\tr T)
	\sum_{\substack{\underline{\lambda}\in \underline{\discel} + \lattice^{2\gengenus}
	\\
	q(\underline{\lambda})=\big(\begin{smallmatrix}
	T & \ast \\
	\ast & T
	\end{smallmatrix}\big)
	}}
	\KMschwartznew{\gengenus}(\underline{\lambda})\nu,
	\]
	and it is the Fourier coefficient of~$\theta_{\underline{\discel}}(\tau_1,\tau_2,\KMschwartznew{\gengenus})\nu$ of index~$(T,T)$ at~$\tau_1=\tau_2=iI_\gengenus$.
	The convergence of~$\lVert\zeta(T,\discel)\rVert_2$ follows from the termwise convergence of the integral over~$X$ of the Fourier coefficients of~$\theta_{\underline{\discel}}(\tau_1,\tau_2,\KMschwartznew{\gengenus})\nu$.
	This is again guaranteed by the Weil convergence criterion.
	
	We now consider the Fourier expansion of~$[\KMtheta{\gengenus}]$ under the assumption that~$\gengenus<(n-1)/2$.
	We know from~\eqref{eq;3isoofcoho} that~$\sqH^{2\gengenus}(X)\cong H^{2\gengenus}(X)$.
	This implies that all square-integrable representatives of~$\drcoho{Z(T,\discel)}$ yield the same~$L^2$-cohomology class.
	The fact that~$[\KMtheta{\gengenus}]$ is a holomorphic Siegel modular form with Fourier expansion~\eqref{eq;FexpKMthetaL2} follows from the analogous properties of~$\drcoho{\KMtheta{\gengenus}}$ proved in~\cite{kudlamillson-intersection}; see Section~\ref{sec;KMthetabasics}.
	\end{proof}
	We define the \emph{volume} of a special cycle~$Z(T,\discel)$ to be
	\[
	\Vol(Z(T,\discel))
	\coloneqq
	\int_{Z(T,\discel)}\omega^{n-\rk(T)}.
	\]
	The integral~$\int_X \omega^{n-\gengenus}\wedge\KMtheta{\gengenus}$ is the \emph{generating series of volumes} of special cycles of codimension~$\genus$ in~$X$.
%
	The following result of Kudla shows that this generating series of volumes is a Siegel Eisenstein series; see~\cite[Theorem~$4.20$ and Corollary~$4.24$]{kudla-speccyc} \textcolor{\newcor}{for a statement in greater generality}.
	Recall the normalized Siegel Eisenstein series~$E^k_{\gengenus,L}$ from Definition~\ref{def;KlingEisseries}.
	\begin{thm}[Kudla]\label{thm;fromKudlagenservol}
	Let~$X$ be an orthogonal Shimura variety arising from an even indefinite lattice~$L$ of signature~$(n,2)$, and let~$k=1+n/2$.
    If~\textcolor{\newcor}{$\gengenus<n/2$},
    then
	\be\label{eq;gengen2inprodkmtheta}
	\int_X \KMtheta{\gengenus}
	\wedge
	\omega^{n-\gengenus}
	=
	(-1)^\gengenus \Vol(X) E^k_{\gengenus,L}.
	\ee
	\end{thm}
	\begin{proof}	
	Let~$\Ab_f$, resp.~$\Ab$, be the ring of finite adeles, resp.~adeles, over~$\QQ$ and let~$V=L\otimes\QQ$.
	For every~$\discel\in D_\lattice^\gengenus$, we may regard the basis vector~$\mathfrak{e}_\discel$ of~$\CC[D_\lattice^\gengenus]$ as the characteristic function~$\varphi_\discel$ of~$\discel + \widehat{L}^\gengenus$ in~$V(\Ab_f)^\gengenus$, where~$\widehat{L}^\gengenus=L^\gengenus\otimes\widehat{\ZZ}$.
	The space of Schwartz--Bruhat functions on~$V(\Ab)^\gengenus$ is denoted by~$\mathcal{S}\big(V(\Ab)^\gengenus\big)$.
	
	We now follow the argument of~\cite[p.~$640$]{bruinier-yang09} and~\cite[pp.~$324$-$325$]{kudla-integrals}, given in genus~$\gengenus=1$, generalizing it to~$\gengenus\ge1$ by means of~\cite{kudla-extensions}.
	For~$s\in\CC$, let~$I_\gengenus(s,\chi_V)$ be the principal series representation of~$\Sp_{2\gengenus}(\Ab)$ induced by~$\chi_V |\cdot |^s$, where~$\chi_V$ is the quadratic character of~$\Ab^*/\QQ^*$ associated to~$V$.
	For every standard section~$\Phi(s)\in I_\gengenus(s,\chi_V)$ we define the Eisenstein series
	\[
	E(g,s;\Phi) \coloneqq \sum_{\gamma\in P(\QQ)\backslash \Sp_{2\gengenus}(\QQ)} \Phi(\gamma g,s),\qquad g\in\Sp_{2\gengenus}(\Ab),
	\]
	where~$P(\QQ)$ is the Siegel parabolic subgroup of~$\Sp_{2\gengenus}(\QQ)$.
	The series~$E(g,s;\Phi)$ converges absolutely to an automorphic form on~$\Sp_{2\gengenus}(\Ab)$ if~$\Re(s)>(\gengenus+1)/2$.
	
	Consider the~$\CC[D_\lattice^\gengenus]$-valued Eisenstein series
	\ba\label{eq;vveisseriesasinkudlabt}
	E(\tau,s;\ell)
	\coloneqq
	\sum_{\discel\in D_\lattice^\genus}\det y^{-\ell/2} E\big(g_\tau,s;\Phi_\infty^\ell \otimes\lambda_f(\varphi_\discel)\big)\mathfrak{e}_\discel,
	\ea
	where $\Phi_\infty^\ell$ is the standard section of the principal series representation of~$\Sp_{2\gengenus}(\RR)$ arising from the character~$\theta\mapsto e^{i\ell\theta}$, and~$\lambda_f\colon\mathcal{S}\big(V(\Ab)^\gengenus\big)\to I_\gengenus(s_0,\chi_V)$ is the intertwining map such that~$\lambda_f(\varphi_\discel)(g_f)=\big(\omega(g_f)\varphi_\discel\big)(0)$ for every~$g_f\in \Sp_{2\gengenus}(\Ab_f)$; see~\cite[Section~$2.1$]{bruinier-yang09} for details.
	
	Let~$s_0=(n-\gengenus+1)/2$.
	As a consequence of the Siegel--Weil formula, Kudla proved in~\cite[Theorem~$4.23$]{kudla-integrals} and~\cite[Theorem~$4.1$]{kudla-speccyc} that the~$\discel$-component of~$E(\tau,s_0;k)$ is a \emph{holomorphic} Eisenstein series and equals the \emph{generating series of volumes}
	\[
	\frac{(-1)^\gengenus}{\Vol(X)}
	\int_X\omega^{n-\gengenus}\wedge(\KMtheta{\gengenus})_\discel,
	\]
	where~$(\KMtheta{\gengenus})_\discel$ is the~$\discel$-component of~$\KMtheta{\gengenus}$.
	
	Following the argument of~\cite[p.~$641$]{bruinier-yang09} and~\cite[Section~IV.$2$]{kudla-extensions}, the vector-valued Eisenstein series~$E(\tau,s_0;k)$ may be rewritten as a Poincaré series arising from the basis vector~$\mathfrak{e}_0$ of~$\CC[D_\lattice^\gengenus]$, namely the Eisenstein series~$\eiszero{k}{\gengenus}{\lattice}$; cf.~\cite[$(2.17)$]{bruinier-yang09}.
	\end{proof}

	Let~$0<\gengenus<(n-1)/2$.
	Recall that~$\sqH^{\gengenus,\gengenus}(X)$ is endowed with a Hermitian inner product~$(\cdot{,}\cdot)$ induced by the Kähler form~$\omega$.
	We may use~\eqref{eq;formmatchonHggsq} and Theorem~\ref{thm;fromKudlagenservol} to compute the orthogonal projection of~$[\KMtheta{\gengenus}]$ to~$\lefschetz^\gengenus \sqH^0(X)=\CC[\omega]^\gengenus$ as
	\[
	\frac{\big([\KMtheta{\gengenus}] , [\omega]^\gengenus\big)}{\big([\omega]^\gengenus,[\omega]^\gengenus\big)}[\omega^\gengenus]
	=
	(-1)^\gengenus E^k_{\gengenus,\lattice}\cdot [\omega]^\gengenus.
	\]
    This implies the following result.
    \begin{cor}\label{rem;projthetafcohoaftervol}
    \textcolor{\newcor}{If~$0<\gengenus<(n-1)/2$, then}
        \be\label{eq;gengen2inprodkmtheta2}
	[\KMtheta{\gengenus}] + (-1)^{\gengenus+1} E^k_{\gengenus,\lattice}\cdot [\omega]^\gengenus \in M^k_{\gengenus,\lattice}\otimes\big(\CC[\omega]^\genus\big)^\perp.
	\ee
    \end{cor}


	\section{The modular Lefschetz decomposition of~$[\KMtheta{\genus}]$}\label{sec;gen2}
	In this section we show that the genus~$\gengenus$ Kudla--Millson theta function~$[\KMtheta{\gengenus}]$ is an almost cusp form, and compute its modular decomposition.
	We also show that this modular decomposition is of geometric interest, since it simultaneously provides the Lefschetz decomposition of~$[\KMtheta{\gengenus}]$ as a cohomology class.
	
	Recall from Section~\ref{sec;mod_dec} that we denote by~$\almcusp{k}{\genus}{\lattice}$ the space of almost cusp forms of weight~$k$ and genus~$\genus$ with respect to the Weil representation~$\weil{\lattice}{\genus}$.
	Its modular decomposition for~$k>2\genus$ is
	\[
	\almcusp{k}{\genus}{\lattice}
	=
	\CC\eiszero{k}{\genus}{L}
	\oplus
	\vklisub{k}{\gengenus}{1}{\lattice}
	\oplus
	\cdots
	\oplus
	\vklisub{k}{\gengenus}{\genus-1}{\lattice}
	\oplus
	S^k_{\genus,\lattice}.
	\]
	where~$\vklisub{k}{\gengenus}{r}{\lattice}$ is the space of vector-valued Klingen Eisenstein series arising from cusp forms in~$S^k_r$.
	Recall also that the Lefschetz decomposition of~$\sqH^{\gengenus,\gengenus}(X)$ is
	\[
	\sqH^{\gengenus,\gengenus}(X)
	=
	\CC[\omega]^\gengenus \oplus \lefschetz^{\gengenus-1} \sqHprim^{1,1}(X) \oplus \cdots \oplus \lefschetz \sqHprim^{\gengenus-1,\gengenus-1}(X)\oplus \sqHprim^{\gengenus,\gengenus}(X).
	\]

	\subsection{Comparison of Siegel and Lefschetz operators}
	Let~$\Phi_\disceltwo\colon M^k_{\gengenus,\lattice}\to M^k_{\gengenus-1,\lattice}$ be the Siegel operator of index~$\disceltwo\in D_\lattice$ introduced in Section~\ref{sec;Siegelphi}.
	We denote the natural extension of that operator to~$M^k_{\gengenus,\lattice}\otimes \sqH^{\gengenus,\gengenus}(X)$ with the same symbol.
	\begin{lemma}\label{lemma;Siegopkmtheta2}
	If~$0<g<(n -1)/2$, then the image of~$[\KMtheta{\gengenus}]$ under the operator~$\Phi_\disceltwo$ is
	\[
	\Phi_\disceltwo\big([\KMtheta{\gengenus}]\big)
	=
	\begin{cases}
	- \lefschetz\big(
	[\KMtheta{\gengenus-1}]
	\big), & \text{if~$\disceltwo=0$,}
	\\
	0, & \text{otherwise.}
	\end{cases}
	\]
	\end{lemma}
	\begin{proof}
	By Lemma~\ref{lemma;L2cohoKMtheta&fexp} the~$L^2$-cohomology class~$[\KMtheta{\gengenus}]$ is a holomorphic Siegel modular form with respect to~$\weil{\lattice}{\genus}$ with Fourier expansion~\eqref{eq;FexpKMthetaL2}.
	By Lemma~\ref{lemma;FexpofimPhiop} we may compute
	\begin{align*}
	\Phi_\disceltwo\big([\KMtheta{\gengenus}]\big)
	&=
	\sum_{\discel\in D_\lattice^{\genus-1}}
	\sum_{\substack{T\in q(\discel)+\halfint{\genus-1} \\ T\ge 0}} \left[ Z\left(\left(\begin{smallmatrix}
	T & 0 \\ 0 & 0
	\end{smallmatrix}\right), (\discel,\disceltwo)\right)\right]\wedge[-\omega]^{\gengenus - \rk(T)} q^{T}\mathfrak{e}_\discel.
	\end{align*}
	It is easy to see from~\eqref{eq;defspcy} that
	\[
	Z\left(\left(\begin{smallmatrix}
	T & 0 \\ 0 & 0
	\end{smallmatrix}\right),(\discel,\disceltwo)\right)
	=
	\begin{cases}
	Z(T,\discel), & \text{if~$\beta=0$,}
	\\
	\emptyset, & \text{otherwise.}
	\end{cases}
	\]
	This implies that if~$\beta\neq 0$, then~$\Phi_\disceltwo\big([\KMtheta{\gengenus}]\big)=0$, and that
	\begin{align*}
	\Phi_0\big([\KMtheta{\gengenus}]\big)
	&=
	- \Big(
	\sum_{\discel\in D_\lattice^{\genus-1}}
	\sum_{\substack{T\in q(\discel)+\halfint{\genus-1} \\ T\ge 0}}
	[Z(T,\discel)]\wedge[-\omega]^{\gengenus -1 - \rk(T)} q^{T}\mathfrak{e}_\discel
	\Big) \wedge [\omega]
	\\
	&= - [\KMtheta{\gengenus-1}]\wedge[\omega]
	= - \lefschetz\big(
	[\KMtheta{\gengenus-1}]
	\big). \qedhere
	\end{align*}
	\end{proof}

	\subsection{The Kudla--Millson lift and its geometric interpretation}
	The \emph{Kudla--Millson lift of genus~$g$} is the linear map defined as
	\be\label{eq;KMliftdef}
	\KMlift{\gengenus}\colon S^k_{g,\lattice}\longrightarrow \sqH^{2\gengenus}(X,\CC),\qquad f\longmapsto\pet{f}{[\KMtheta{\gengenus}]}.
	\ee
	It has been introduced by Kudla and Millson in~\cite{kudlamillson-tubes}, and may be considered as a bridge to study the geometric properties of~$X$ in terms of the arithmetic properties of Siegel modular forms.
	
	The following result, essential for the proof of the Lefschetz decomposition of~$[\KMtheta{\genus}]$, provides a geometric interpretation of the image of~$\KMlift{\genus}$.
	\begin{thm}\label{thm;imageofKMlift}
	If~$\genus < (n+2)/4$, then	the image of the Kudla--Millson lift~$\KMlift{\gengenus}$ is contained in the primitive cohomology group~$\sqHprim^{\gengenus,\gengenus}(X)$.
	\end{thm}
	\begin{proof}
	We begin with the case of~$\Gamma=\disk{}^+(L)$.
	Let~$K$ be the open compact subgroup of~$G(\Ab_f)$ as in Section~\ref{sec;adelicspcy}.
	By Theorem~\ref{thm;noetherlef}, the space~$\sqH^{r,r}(X)$ is generated by the cohomology classes of adelic special cycles for every~$r\leq\genus-1$.
	Recall that for every~$\varphi\in\mathcal{S}\big(V^{r}(\Ab_f)\big)$, there exists an even lattice~$M\subset V$ such that~$\varphi\in\mathcal{S}_{M^{r}}$; see e.g.~\cite{kudla-speccyc}.
	Hence, for every~$r\leq\genus-1$ there exists a sufficiently small even lattice~$M_r\subset V$ such that~$\sqH^{r,r}(X)$ is generated by the adelic special cycles of the form
	\[
	[Z(T,\varphi,K)]\wedge[\omega^{r-\rk T}],
	\qquad
	\text{where~$T\in\Sym_r(\QQ)$ and~$\varphi\in(\mathcal{S}_{(M_r)^r})^K$.}
	\]
	Without loss of generality, we may assume that the~$M_r$ are all equal to the same even lattice~$M\subset V$, and also that~$M\subset L$.
	
	Saying that the image of~$\KMlift{\genus}$ is contained in~$\sqHprim^{\genus,\genus}(X)=\big(\lefschetz \sqH^{\genus-1,\genus-1}(X)\big)^\perp$ is equivalent to saying that such an image is orthogonal to~$\lefschetz^{\genus-\rk T}\big([Z(T,\varphi,K)]\big)$ for every~$T\in \Sym_{\genus-1}(\QQ)$ and~$\varphi\in (\mathcal{S}_{M^{\genus-1}})^K$.
	
	Let~$X_M=\disk{}^+(M)\backslash\domain$, and let~$\pi\colon X_M\to X$ be the quotient map.
	We denote by~$\KMtheta{\genus-1}^M$ the Kudla--Millson theta function arising from the lattice~$M$.
	Let~$f\in S^k_{\genus,\lattice}$.
	It is easy to check by using~\eqref{eq;formmatchonHggsq} and induction on~$\genus$ that if
	\be\label{eq;wwwtpingengenuspldu}
	\int_{X_M} \pi^*\big(\KMlift{\genus}(f)\big)\wedge \KMtheta{\genus-1}^M \wedge \omega^{n-(2\genus-1)}
	=0,
	\ee
	then
	\be\label{eq;wwwtpingengenusplduimpl}
	\big(
	\KMlift{\genus}(f)
	,
	\lefschetz^{\genus-\rk T}[Z(T,\varphi,K)]
	\big)=0
	\qquad
	\text{for every~$T$ and~$\varphi$ as above}.
	\ee
	To show that the image of~$\KMlift{\genus}$ lies in~$\sqHprim^{\genus,\genus}(X)$, it is then enough to prove~\eqref{eq;wwwtpingengenuspldu}.
	
	By Lemma~\ref{lemma;genofBYandS} and Remark~\ref{rem;resandpullbKMthetafunct} we deduce that
	\[
	\pi^*\big(\KMlift{\genus}(f)\big)
	=
	\bigpet{f}{\pi^*([\KMtheta{\genus}])}
	=
	\bigpet{f}{\tr_{\genus,\lattice/M}([\KMtheta{\genus}^M])}
	=
	\bigpet{f_M}{[\KMtheta{\genus}^M]}.
	\]
	By the Fubini--Tonelli Theorem, we may then rewrite the integral in~\eqref{eq;wwwtpingengenuspldu}, keeping track of the variable~$\tau_1\in\HH_{\genus-1}$, as
	\begin{align*}
	&\int_{X_M}
	\pi^*\big(\KMlift{\genus}(f)\big)
	\wedge
	\KMtheta{\genus-1}^M(\tau_1)
	\wedge
	\omega^{n-(2\genus-1)}
	\\
	&\quad
	=
	\Bigpet{f_M}{\int_{X_M} \KMtheta{\genus-1}^M(-\overline{\tau_1}) \wedge \KMtheta{\genus}^M \wedge \omega^{n-(2\genus-1)}},
	\end{align*}
	where to deduce the equality above we used that~$\KMtheta{\genus-1}$ is a modular form with \emph{real} Fourier coefficients, and hence~$\overline{\KMtheta{\genus-1}(\tau_1)} = \KMtheta{\genus-1}(-\overline{\tau_1})$.
	By~\eqref{eq;wedgeprodformforKMtheta} and Theorem~\ref{thm;fromKudlagenservol} we may continue the computation above as
	\bas
	\,&\int_{X_M} \pi^*\big(\KMlift{\genus}(f)\big)\wedge
	\KMtheta{\genus-1}^M(\tau_1)\wedge
	\omega^{n-(2\genus-1)}
	\\
	&
	\quad
	=
	\Bigpet{f_M}{\int_{X_M} \KMtheta{2\genus-1}^M\big(\begin{smallmatrix}
	-\overline{\tau_1} & \\
	 & \ast
	\end{smallmatrix}\big) \wedge \omega^{n-(2\genus-1)}}
	\\
	&\quad
	=
	(-1)^{2\genus-1}
	\Vol(X_M)
	\Bigpet{f_M}{E^k_{2\genus-1,M}\big(\begin{smallmatrix}
	-\overline{\tau_1} & \\
	 & \ast
	\end{smallmatrix}\big)},
	\eas
	where we denoted by~$\ast$ the variable of integration in~$\HH_{\genus}$ for the Petersson inner product.
	
	By Corollary~\ref{cor;vanishingGmapvv} we deduce that
	\[
	\Bigpet{f_M}{E^k_{2\genus-1,M}\big(\begin{smallmatrix}
	-\overline{\tau_1} & \\
	 & \ast
	\end{smallmatrix}\big)} = 0,
	\]
	which then implies that~$\KMlift{\genus}(S^k_\genus)\subseteq \sqHprim^{\gengenus,\gengenus}(X)$.	
	
	We now consider the case where~$\Gamma$ is a subgroup of finite index in~$\Gamma_L\coloneqq\disk{}^+(\lattice)$.
	Let
	\[
	\pi\colon X_\Gamma\coloneqq \Gamma\backslash\domain
	\longrightarrow
	X_L\coloneqq\Gamma_\lattice\backslash\domain
	\]
	be the projection map.
	The pullback~$\pi^*\colon H^{2\genus}(X_\lattice)\to H^{2\genus}(X_\Gamma)\cong \sqH^{2\genus}(X_\Gamma)$ is injective and identifies~$H^{2\genus}(X_\lattice)$ with the subspace~$\sqH^{2\genus}(X_\Gamma)^{\Gamma\backslash\Gamma_\lattice}$ of~$\Gamma\backslash\Gamma_\lattice$-invariant classes in~$\sqH^{2\genus}(X_\Gamma)$.
	The cohomology class of the Kudla--Millson theta function attached to~$L$ induced on~$X_\Gamma$ is the pullback~$\pi^*\big([\KMtheta{\genus}]\big)\in \sqH^{2\genus}(X_\Gamma)^{\Gamma\backslash\Gamma_\lattice}$ of the one defined on~$X_\lattice$.
	Since the Kähler form~$\omega$, and hence the Hermitian inner product on~$\sqH^{2\genus}(X_\Gamma)$, is~$\Gamma_\lattice$-invariant, we deduce that
	\[
	\sqH^{2\genus}(X_\Gamma)^{\Gamma\backslash\Gamma_\lattice}
	=
	\bigoplus_{p+q=2\genus} \sqH^{p,q}(X_\Gamma)^{\Gamma\backslash\Gamma_\lattice}
	\qquad
	\text{and}
	\qquad
	\sqH^{\genus,\genus}(X_\Gamma)^{\Gamma\backslash\Gamma_\lattice}
	=
	\bigoplus_{i\geq 0}
	\lefschetzî\sqHprim^{r-2i}(X_\Gamma)^{\Gamma\backslash\Gamma_\lattice}.
	\]
Since the pullback~$\pi^*$ maps primitive cohomology groups of~$X_L$ into primitive cohomology groups of~$X_\Gamma$, the image of~$\KMlift{\genus}$ lies in primitive cohomology also on~$X_\Gamma$.
	\end{proof}	
	
	\subsection{The modular Lefschetz decomposition in genus 1}\label{sec;gen1}
	In this section we describe the modular and Lefschetz decompositions of the class~$[\KMtheta{1}]$ of the genus~$1$ Kudla--Millson theta function.
	
%
	Let~$n\geq 3$.
	By~\eqref{eq;3isoofcoho}, as well as~\cite[Korollar~$4.2$]{weissauer-88} for the special case of~$n=3$, we know that~$\sqH^{1,1}(X)\cong H^{1,1}(X)$.
	By Corollary~\ref{cor;blmminnonadcase} and~\cite{HeHoffman}, the space~$\sqH^{1,1}(X)$ is generated by the cohomology classes of special divisors.
	The class~$[\KMtheta{1}]$ is the generating series of cohomology classes of such divisors, more precisely
	\bes
	[\KMtheta{1}]=-[\omega]\mathfrak{e}_0 + \sum_{\discel\in D_\lattice}\sum_{\substack{m\in q(\discel)+\ZZ \\ m>0}} [Z(m,\discel)]q^m \mathfrak{e}_\discel.
	\ees
	
	The \emph{Lefschetz decomposition} of~$[\KMtheta{1}]$ is a trivial consequence of Corollary~\ref{rem;projthetafcohoaftervol}, since~\eqref{eq;gengen2inprodkmtheta2} implies for~$g=1$ that
	\be\label{eq;lefdecgen1coho}
	[\KMtheta{1}]
	=
	\Big(\underbrace{- E^k_{1,\lattice} \cdot [\omega]}_{\text{in }M^k_{1,\lattice}\otimes{}\CC[\omega]}\Big)
	+
	\Big(\underbrace{[\KMtheta{1}] + E^k_{1,\lattice} \cdot [\omega]}_{\text{in }M^k_{1,\lattice}\otimes \sqHprim^{1,1}(X)}\Big).
	\ee
	
	Since the Fourier coefficient of index~$(0,\discel)$ of~$[\KMtheta{1}] + E^k_{1,\lattice} \cdot [\omega]$ is trivial for every~$\discel\in D_\lattice$, the primitive part of~$[\KMtheta{1}]$ lies in~$S^k_{1,\lattice}\otimes \sqHprim^{1,1}(X)$, i.e.~it behaves as a cusp form with respect to the variable~$\tau\in\HH$.
	The part in~$\CC[\omega]$ of~$[\KMtheta{1}]$ lies in~$\CC E^k_{1,\lattice} \otimes \CC[\omega]$, namely it behaves as an Eisenstein series with respect to the variable~$\tau$.
	This implies that the Lefschetz decomposition of~$[\KMtheta{1}]$ with respect to the orthogonal variable is also the modular decomposition of~$[\KMtheta{1}]$ with respect to the symplectic variable.

	\subsection{The case of general genus}\label{sec;cogengenus}
	We are now ready to prove the main result of this article.
	Recall from~\eqref{eq;defofKlinblamap} that the map~$\vklin{k}{\genus}{\genus-1}{\lattice}$ is the inverse of the restriction of the Siegel operator~$\Phi_0$ to~$\vklisub{k}{\gengenus}{0}{\lattice}\oplus\cdots\oplus \vklisub{k}{\gengenus}{\genus-1}{\lattice}$.
	\begin{thm}\label{thm;lefschetzdecgengenus}
	Let~$\gengenus<(n+2)/4
	$ be a positive integer.
	Then~$[\KMtheta{\gengenus}]$ is an almost cusp form.
	Its primitive part is
	\be\label{eq;primpartgengenus}
	[\KMtheta{\gengenus}]_0 = [\KMtheta{\gengenus}]
	+
	\lefschetz\circ\vklin{k}{\genus}{\genus-1}{\lattice}\big([\KMtheta{\gengenus-1}]\big),
	\ee
	and lies in~$S^k_{\gengenus,\lattice}\otimes \sqHprim^{\gengenus,\gengenus}(X)$.
	In particular, it is a Siegel cusp form with respect to the symplectic variable.
	The \emph{Lefschetz decomposition} of~$[\KMtheta{\gengenus}]$ equals the \emph{modular decomposition} of~$[\KMtheta{\gengenus}]$, and is given by
	\ba\label{eq;lefschetzdecgengenus}
	\,[\KMtheta{\gengenus}]
	=
	\sum_{r=0}^{\gengenus} (-\lefschetz)^{\gengenus-r}\circ
	\klinser{k}{\gengenus}{r}{\lattice} \big([\KMtheta{r}]_0\big).
	\ea
	Here the summand of index~$r$
	lies in~$\vklisub{k}{\gengenus}{r}{\lattice} \otimes \lefschetz^{\gengenus-r} \sqHprim^{r,r}(X)$, in particular it is a Klingen Eisenstein series arising from a genus~$r$ Siegel cusp form.
	
	All previous assertions hold also if we replace~$L^2$-classes by de Rham cohomology classes.
	\end{thm}
	\begin{proof}[Proof of Theorem~\ref{thm;lefschetzdecgengenus}]
	We prove the result by induction on the genus~$\gengenus$.
	The case of genus~$1$ has already been illustrated in Section~\ref{sec;gen1}.
	Suppose that~$\genus>1$ and that Theorem~\ref{thm;lefschetzdecgengenus} is valid for every genus less than~$g$.
	We prove that it is valid also for~$\gengenus$.

	To prove that the right-hand side of~\eqref{eq;primpartgengenus} is a cusp form, it is enough to show that it lies in the kernel of~$\Phi_\disceltwo\colon M^k_{g,\lattice}\otimes \sqH^{g,g}(X)\to M^k_{g-1,\lattice}\otimes \sqH^{g,g}(X)$ for every~$\disceltwo\in D_\lattice$; see Lemma~\ref{lemma;FexpofimPhiop}.
	This follows easily from Lemma~\ref{lemma;Siegopkmtheta2} and Theorem~\ref{thm;Phi0onkling}.
	
	It is easy to check that the right-hand side of~\eqref{eq;lefschetzdecgengenus} equals~$[\KMtheta{\genus}]$ under the assumption of~\eqref{eq;primpartgengenus}.

	The inductive hypothesis that~$[\KMtheta{r}]_0$ lies in~$S^k_{r,L}\otimes \sqHprim^{r,r}(X)$ for every~$r=0,\dots,\gengenus-1$ implies that~$\lefschetz^{\genus-r}\circ\klinser{k}{\gengenus}{r}{\lattice}\big([\KMtheta{r}]_0\big)$ lies in~$\vklisub{k}{\gengenus}{r}{\lattice}\otimes \lefschetz^{\gengenus-r} \sqHprim^{r,r}(X)$.
	It remains to show~\eqref{eq;primpartgengenus}, assuming by induction that~\eqref{eq;primpartgengenus} is satisfied for every genus less than~$g$.
	Since the right-hand side of~\eqref{eq;primpartgengenus} is the cuspidal part of~$[\KMtheta{\gengenus}]$, it is enough to show that the Kudla--Millson lift
	\[
	\KMlift{\genus}\colon S^k_{\genus,\lattice}\longrightarrow \sqH^{\genus,\genus}(X),\qquad f\longmapsto\pet{f}{[\KMtheta{\genus}]}
	\]
	has image contained in~$\sqHprim^{\genus,\genus}(X)$.
	In fact, the space~$S^k_{\genus,\lattice}$ is generated by Poincaré series, hence we may deduce from Lemma~\ref{lemma;petSiegpoinandcusp} that every Fourier coefficient of the right-hand side of~\eqref{eq;primpartgengenus} is contained in~$\sqHprim^{\genus,\genus}(X)$.
	The primitivity of the image of~$\KMlift{\genus}$ follows from Theorem~\ref{thm;imageofKMlift}.
	\end{proof}

	\section{Geometric applications of the Lefschetz decomposition of~$[\KMtheta{\genus}]$}\label{sec;appgen1}
	In this section we illustrate some geometric applications of the Lefschetz decomposition of~$[\KMtheta{\genus}]$ and the geometric interpretation of the lift~$\KMlift{\genus}$.
	
	\subsection{More on the image of the Kudla--Millson lift}\label{sec;imageKMlift}
	From the Lefschetz decomposition of~$[\KMtheta{\gengenus}]$ we may sharpen Theorem~\ref{thm;imageofKMlift} as follows.
	\begin{cor}\label{cor;imageofKMliftissurj}
	Let~$\genus < (n+2)/4$.
	If the codimension~$\genus$ special cycles on~$X$ arising from~$L$ span the whole~$\sqH^{\genus,\genus}(X)$, then the image of~$\KMlift{\gengenus}$ is~$\sqHprim^{\gengenus,\gengenus}(X)$.
	\end{cor}
	A criterion on~$L$ for which~$H^{\genus,\genus}(X)$ is spanned by the codimension~$\genus$ special cycles is given by Corollary~\ref{cor;blmminnonadcase} for small~$g$.
	\begin{proof}
%
	By Lemma~\ref{lemma;petSiegpoinandcusp} the space~$S^k_{\gengenus,\lattice}$ is generated by Poincaré series~$\poincv{\gengenus}{\discel}{T}{\lattice}$, where~$\discel\in D_\lattice^\genus$ and~$T\in q(\discel)+\halfint{\genus}$ is positive definite.
	Together with Lemma~\ref{lemma;petSiegpoinandcusp} and Theorem~\ref{thm;lefschetzdecgengenus}, this implies that the image of~$\KMlift{\genus}$ is generated by the primitive classes~$[Z(T,\discel)]_0$ of codimension~$\genus$ special cycles, with~$T>0$.
	These primitive classes span the whole~$\sqHprim^{\gengenus,\gengenus}(X)$ by assumption.
%
	\end{proof}
	
	In~\cite[Theorem, p.~$4$]{kudlamillson-tubes} Kudla and Millson claim that the image of~$\KMlift{\gengenus}$ is the span in~$\sqH^{\gengenus,\gengenus}(X)$ of the classes of codimension~$\gengenus$ special cycles.
	This is in contrast with Theorem~\ref{thm;imageofKMlift} and Corollary~\ref{cor;imageofKMliftissurj}, since that span does not lie in~$\sqHprim^{\gengenus,\gengenus}(X)$.
	In fact, the proof of~\cite[Theorem, p.~$4$]{kudlamillson-tubes} contains an error.
	The goal of the remaining part of this section is to illustrate the idea behind the proof of~\cite[Theorem, p.~$4$]{kudlamillson-tubes} and to provide a counterexample to it.
	The latter does not rely on the main results of this article. 
	\\
	
	To prove~\cite[Theorem, p.~$4$]{kudlamillson-tubes} Kudla and Millson used the following incorrect claim, which for simplicity we state only in the case of genus~$g=1$ and under the assumption that~$L$ is \emph{unimodular}:
	For every closed compactly supported~$(n-1,n-1)$-differential form~$\eta$ on~$X$, if~$\pet{\theta_\eta}{P_m}=0$, then~$\int_{Z(m)}\eta=0$.
	Here~$\theta_\eta\coloneqq\int_X \eta\wedge\KMtheta{1}\in M^k_1$,~$P_m$ is the normalized Poincaré series of positive index~$m$, and~$Z(m)$ is the~$m$-th special divisor on~$X$.
	
	Their attempt of proof of the claim above is based on an imprecise analytic continuation, which we now explain.
	Define the \emph{non-holomorphic Poincaré series of index~$m$} as
	\[
	P_m(\tau,s)
	= 
	c(s)
	\sum_{\gamma\in\Gamma_\infty\backslash\SL_2(\ZZ)}
	(y^sq^m)|_k\gamma
	,\qquad s\in \CC,
	\]
	where~$c(s)$ is the normalization constant such that~$\pet{f}{P_m(\cdot,s)}=c_m(f)$ for every~${f\in S^k_1}$, and~$\Gamma_\infty$ is the parabolic subgroup of translations in~$\SL_2(\ZZ)$.
	Note that~$P_m(\tau,0)=P_m(\tau)$.

	As illustrated in~\cite[p.~$25$]{kudlamillson-tubes}, if~$\Re(s)>1$, then one can apply the unfolding trick to deduce that
	\begin{align}\label{eq;unfoldingKMmistend}
	\begin{split}
	\bigpet{\theta_\eta}{P_m(\cdot,s)}
	=
	\frac{c(s)\cdot\Gamma(k+s-1)}{(4\pi m)^{k+s-1}}\int_{Z(m)}\eta.
	\end{split}
	\end{align}
	Although~\eqref{eq;unfoldingKMmistend} is justified only for~$\Re(s)>1$, we may consider the analytic continuation of the right-hand side of~\eqref{eq;unfoldingKMmistend} and evaluate it at~$s=0$.
	This equals~$\int_{Z(m)}\eta$.
	
	The problem is that such value does not equal~$\bigpet{\theta_\eta}{P_m(\cdot,0)}$, and in fact the vanishing of~$\bigpet{\theta_\eta}{P_m(\cdot,0)}$ does not imply the vanishing of~$\int_{Z(m)}\eta$.
	This is shown in the following result, based on the fact that~$\bigpet{\theta_\eta}{P_m(\cdot,0)}=c_m(\theta_\eta)$ only if~$\theta_\eta$ is a cusp form, which is not the case for certain~$\eta$.
	\begin{lemma}
	The vanishing of~$\bigpet{\theta_\eta}{P_m(\cdot,0)}$ does not imply the vanishing of~$\int_{Z(m)}\eta$.
	\end{lemma}
	\begin{proof}
%
	We rewrite the weight~$k$ modular form~$\theta_\eta$ with respect to the modular decomposition~${M^k_1=\CC E^k_1\oplus S^k_1}$ in Eisenstein and cuspidal parts as~$\theta_\eta=a_\eta E^k_1 + \tilde{\theta}_\eta$, for some constant~$a_\eta\in\CC$ and some cusp form~$\tilde{\theta}_\eta\in S^k_1$.
	
	Let~$\eta$ be such that~$a_\eta\neq 0$.
	An example of such a~$\eta$ is given by the Poincaré dual of a compact orthogonal Shimura curve, as follows.
	Let~$X=\Gamma\backslash \domain$, where~$\domain$ is the Hermitian symmetric domain associated to some unimodular lattice~$L$ of signature~$(n,2)$.
	The orthogonal Shimura subvarieties of codimension~$r$ in~$X$ arise from Hermitian symmetric subdomains~$\domain_W$ in~$\domain$ defined as the orthogonal complement in~$\domain$ of some positive definite~$r$-dimensional subspace~$W\subseteq L\otimes\QQ$ .
	We choose~$r=n-1$ and~$W$ such that~$W^\perp$ is anisotropic in~$L\otimes\QQ$.
	The subvariety arising from~$\domain_W$ is then a \emph{compact} curve~$C$ in~$X$.
	The Poincaré dual of~$C$ is a~$(n-1,n-1)$-differential form~$\eta$ of compact support, and the associated modular form~$\theta_\eta$ has Fourier expansion
	\[
	\theta_\eta(\tau)=-\Vol(C) + \sum_{m>0}\Big(\int_{Z(m)}\eta\Big)q^m.
	\]
	This is not a cusp form, since the constant Fourier coefficient is non-zero.
	
	The left-hand side of~\eqref{eq;unfoldingKMmistend} evaluated at~$s=0$ may be computed as
	\[
	\bigpet{\theta_\eta}{P_m(\cdot,0)}=\bigpet{\tilde{\theta}_\eta}{P_m(\cdot,0)}=c_m(\tilde{\theta}_\eta),
	\]
	since the Eisenstein part of~$\theta_\eta$ is orthogonal to the holomorphic Poincaré series of index~$m$.
	
	The right-hand side of~\eqref{eq;unfoldingKMmistend} evaluated at~$s=0$ is the~$m$-th Fourier coefficient of~$\theta_\eta$.
	This differs from the~$m$-th Fourier coefficient of the cuspidal part~$\tilde{\theta}_\eta$, since~$c_m(E^k_1)\neq 0$ for every~$m$.
%
%
	\end{proof}

	\subsection{Dimensions of cohomology groups}\label{sec;dimcoho}
	Let~$\genus < (n+2)/4$, and let
	\be\label{eq;tocohomap}
	\modcoho{\genus}{\lattice}\colon \almcusp{k}{\genus}{\lattice}\longrightarrow \sqH^{\genus,\genus}(X)
	\ee
	be the linear map constructed on each term of the modular decomposition of~$\almcusp{k}{\genus}{\lattice}$ as
	\[
	\klinser{k}{\gengenus}{r}{\lattice}(f)
	\longmapsto
	(-\lefschetz)^{\genus-r} \circ \KMlift{r}(f)\qquad \text{for all~$f\in S^k_{r,\lattice}$}.
	\]
	As immediate consequence of Theorem~\ref{thm;imageofKMlift}, we may deduce the following result.
	\begin{cor}\label{cor;ifinjthenformcoho}
	Let~$\genus<(n+2)/4$.
	Then
	\[
	\dim\lefschetz^{\genus-r} H_0^{r,r}(X)\geq\dim \KMlift{r}\big(S^k_{r,L}\big)\qquad\text{for every~$0\leq r\leq\genus$}.
	\]
	If the Kudla--Millson lift~$\KMlift{r}$ is injective for every~$r\leq \genus$, then so is~$\modcoho{\genus}{\lattice}$ and
	\[
	\dim\lefschetz^{\genus-r} H_0^{r,r}(X)\geq\dim S^k_{r,\lattice}
	\qquad\text{and}\qquad
	\dim H^{\genus,\genus}(X)\geq\dim \almcusp{k}{\genus}{\lattice}.
	\]
	\end{cor}
%
	\begin{proof}
	By~\eqref{eq;3isoofcoho} and~\cite{weissauer-88} the group~$\sqH^{\genus,\genus}(X)$ is isomorphic to and induces a Lefschetz decomposition on~$H^{\genus,\genus}(X)$.
	The operators~$\lefschetz^{\genus-r}\colon \sqH^{r,r}(X)\to \sqH^{\genus,\genus}(X)$ are injective by Theorem~\ref{thm;lefdecL2}.
	Corollary~\ref{cor;ifinjthenformcoho} is then a consequence of Theorem~\ref{thm;imageofKMlift}.
	\end{proof}
	\begin{rem}
	The dual space~$(\almcusp{k}{\genus}{\lattice})^*$ is generated by the \emph{coefficient extraction functionals}
	\[
	c_{T,\discel}
	\colon \almcusp{k}{\genus}{\lattice}\longrightarrow\CC,
	\qquad
	f\longmapsto c_T(f_\discel).
	\]
	We identify~$\vklisub{k}{\gengenus}{r}{\lattice}$ with its dual space by means of the Petersson product of~$S^k_{r,\lattice}$ as
	\[
	\vklisub{k}{\gengenus}{r}{\lattice}\longrightarrow
	\big(\vklisub{k}{\gengenus}{r}{\lattice}\big)^*,
	\qquad
	\klinser{k}{\gengenus}{r}{\lattice}(f)\longmapsto
	\Big(
	\klinser{k}{\gengenus}{r}{\lattice}(f')\mapsto \pet{f'}{f}
	\Big).
	\]
	It is easy to use such identifications to show that surjectivity and injectivity results on~$\modcoho{\genus}{\lattice}$ imply analogous results on the linear map
	\[
	\dumodcoho{\genus}{\lattice} \colon
	(\almcusp{k}{\genus}{\lattice})^* \longrightarrow \sqH^{\genus,\genus}(X),
	\qquad
	c_{T,\discel}\longmapsto [Z(T,\discel)].
	\]
	Such a map has been used in~\cite{bruiniermöller}, \cite{zuffetti-equid} and~\cite{zuffetti-conescod2} to study cones generated by special divisors in~$\Pic(X)\otimes\RR$ and cones of codimension~$2$ special cycles in~$\CH^{2}(X)\otimes\RR$.
	Important for that purpose is the injectivity of~$\dumodcoho{\genus}{\lattice}$, which has been studied only for~$\genus=1$ with a different approach; see~\cite{bruinierbook} and~\cite{borcherds-grassmannians}.
	By Corollary~\ref{cor;ifinjthenformcoho} we may deduce the injectivity of~$\dumodcoho{\genus}{\lattice}$ from the injectivity of Kudla--Millson lifts.
	This solves the problem raised in~\cite[Remark~$4.13$]{zuffetti-conescod2}.
	\end{rem}
	
	The injectivity of the Kudla--Millson lift of genus~$g=1$ has been studied in several articles.
	The first author proved in~\cite{bruinierbook} and~\cite{bruinier-converse} that~$\KMlift{1}$ is injective whenever the lattice~$L$ splits off a hyperbolic plane and a scaled hyperbolic plane.
	Further works of the first author and Funke~\cite{bruinierfunke}, Stein~\cite{stein-converse} and Metzler and the second author~\cite{metzlerzuffetti} provide further generalizations to locally symmetric spaces of orthogonal type.
	
	The second author illustrated in~\cite{zuffetti-unfolding} a new method to prove the injectivity in genus~$1$, with the advantage of paving the ground for a strategy to deduce injectivity in higher genus.	
	This has been implemented in genus~$2$ by Kiefer and the second author in~\cite{kieferzuffetti}.
	They proved the injectivity of~$\KMlift{2}$ under the assumption that~$\lattice$ splits off two hyperbolic planes, and that the orthogonal complement~$\lattice^+$ of such hyperbolic planes is such that~$\lattice^+\otimes\ZZ_p$ splits off two hyperbolic planes over~$\ZZ_p$ for every prime~$p$.
	An example of~$L^+$ satisfying the latter condition is the~$8$th root lattice~$E_8$; see e.g.~\cite[Satz~$15.8$]{kneser} and~\cite[Theorem~$5.2.2$]{kitaoka}.
	
	The following results follow directly from the injectivity results stated above and Corollaries~\ref{cor;blmminnonadcase}, \ref{cor;imageofKMliftissurj} and~\ref{cor;ifinjthenformcoho}.
	\begin{cor}\label{cor;dimcoho1}
	Let~$X=\disk{}^+(L)\backslash\domain$ with~$L$ of rank~$rk(L)>4$.
	If~$L$ splits off~$U^{\oplus 2}$, then
	\[
	\dim H^{1,1}_0(X)=\dim S^k_{1,\lattice}
	\qquad\text{and}\qquad
	\dim H^{1,1}(X)= \dim \almcusp{k}{1}{\lattice}= 1 + \dim S^k_{1,\lattice}.
	\]
	\end{cor}
	\begin{cor}\label{cor;dimcoho2}
	Let~$X=\disk{}^+(L)\backslash\domain$.
	If~$\lattice$ splits off~$E_8\oplus U^{\oplus 2}$,
	then~$\dim H^{2,2}_0(X)=\dim S^k_{2,\lattice}$, $\dim \lefschetz H^{1,1}_0(X)=\dim S^k_{1,\lattice}$ and
	\[
	\dim H^{2,2}(X)= \dim \almcusp{k}{2}{\lattice}= 1 + \dim S^k_{1,\lattice} + \dim S^k_{2,\lattice}.
	\]
	\end{cor}
	
	\begin{ex}
	Interesting examples of orthogonal Shimura varieties are the (coarse) moduli spaces of irreducible holomorphic symplectic manifolds.
	Several of them satisfy the hypothesis of Corollaries~\ref{cor;dimcoho1} and~\ref{cor;dimcoho2}, for instance the moduli spaces of quasi-polarized K3 surfaces, which are orthogonal Shimura varieties arising from the lattices
	\[
	L_d=\langle 2d \rangle\oplus U^{\oplus 2}\oplus E_8^{\oplus2},
	\]
	where~$d\in\ZZ_{>0}$ is the degree of the polarization.
	Note that we choose the bilinear form on the middle cohomology groups of~K3 surfaces to be the negative of the intersection paring.
	
	For K3 surfaces the degree is the only invariant of a polarization.
	This is not the case for the Hilbert schemes~$S^{[m]}$, where~$S$ is a K3 surface.
	We refer to~\cite[Section~$3$]{GHS-compositio} for a general overview, and restrict here to the \emph{split case}, where the moduli space of deformation~K3${}^{[m]}$ manifolds arises from a lattice of the form
	\[
	L_{t,d}=\langle 2t \rangle\oplus\langle 2d \rangle\oplus U^{\oplus 2}\oplus E_8^{\oplus2}
	\] 
	for some~$t,d\in\ZZ_{>0}$.
	Also in this case the hypothesis of Corollary~\ref{cor;dimcoho2} are satisfied.
	\end{ex}
	
	We remark that explicit formulas to compute~$\dim S^k_{\genus,\lattice}$ are available in the literature only if either~$\genus=1$, or~$L$ is unimodular and~$\genus$ is small; see~\cite[Section~$7$]{borcherds-reflection}, \cite{fischer}, \cite[Section~$5.4$]{taibi}, \cite{ChenevierLannes}.
	The computation of~$\dim S^k_{\genus,\lattice}$ for a general lattice~$L$ is still an interesting open problem.
	
	We remark that the dimension formulas stated in this section hold also for the space spanned by the \emph{rational classes} of special cycles in the Chow group~$\CH^\genus(X)_\CC\coloneqq\CH^\genus(X)\otimes\CC$.
	We define that subspace as
	\[
	\SC^\genus(X)
	\coloneqq
	\Big\langle
	\{Z(T,\discel)\}\cdot\{\omega\}^{\genus-\rk T}
	:
	\text{$\discel\in D_\lattice^\genus$, $T\in q(\discel) + \halfint{\genus}$ and $T\geq 0$}
	\Big\rangle_\CC
	\subseteq\CH^\genus(X)_\CC.
	\]
	As a consequence of Kudla's Modularity Conjecture, proved by the first author and Raum in~\cite{bruinierraum}, it is easy to deduce the following result with the same wording of Corollary~\ref{cor;ifinjthenformcoho}.
	\begin{cor}\label{cor;dimspanCH}
	Let~$X=\disk{}^+(L)\backslash\domain$ and let~$\genus<(n+2)/4$.
	If the Kudla--Millson lift~$\KMlift{r}$ is \emph{injective} for every~$r\leq \genus$, then~$\dim SC^\genus(X)=\dim \almcusp{k}{\genus}{\lattice}$.
	\end{cor}

	\printbibliography

\end{document}